%%%%%%%%%%%%%%%%%%%%%%%%%%%%%%%%%%%%%%%%%%%%%%%%%%%%%%%%%%%%%%%%%%%%%%%%%%%%%%%
%% THIS IS forum.template, A TEMPLATE FOR TeX SOURCE FILES FOR SEMIGROUP FORUM.
%% AUTHOR: OPEN A FILE FOR YOUR TeX SOURCE FILE AND LOAD
%% THIS AT THE BEGINNING OF THAT FILE.  FILL IN THE INDICATED INFORMATION.
%%%%%%%%%%%%%%%%%%%%%%%%%%%%%%%%%%%%%%%%%%%%%%%%%%%%%%%%%%%%%%%%%%%%%%%%%%%%%%%
%% INSERT YOUR OWN MACROS HERE

\def\cm{com\-mu\-ta\-ti\-ve mo\-no\-id}
\def\poag{par\-tial\-ly or\-de\-red a\-be\-lian group}

\def\uip#1{\lceil#1\rceil}
\def\and{{\rm and}}
\def\or{{\rm or}}

\font\sevensl=cmsl10 at 7pt
\font\sevenit=cmti7
\font\bfi=cmbxti10

\font\EUFMX=eufm10 at 10pt
\font\EUFMVII=eufm7 at 7pt 
\font\EUFMV=eufm5 at 5pt 
\textfont8=\EUFMX
\scriptfont8=\EUFMVII\scriptscriptfont8=\EUFMV
\def\goth{\fam8\EUFMX}

\font\MSBMX=msbm10 at 10pt
\font\MSBMVII=msbm7 at 7pt
\font\MSBMV=msbm5 at 5pt 
\textfont9=\MSBMX
\scriptfont9=\MSBMVII\scriptscriptfont9=\MSBMV
\def\bb{\fam9\MSBMX}

\def\berg{1}
\def\clpr{2}
\def\dobb{3}
\def\good{4}
\def\goodII{5}
\def\goodIII{6}
\def\gril{7}
\def\howi{8}

\def\pard{10}
\def\wehr{11}
\def\wehrII{12}
\def\wehrIII{13}
\def\wehrIV{14}

\def\tvi{\vrule height 12pt depth 5pt width 0pt}
\def\tv{\tvi\vrule}
\def\cc#1{\hfill\kern.7em#1\kern.7em\hfill}

\def\rfl#1#2{\smash{\mathop{\hbox to 12mm{\rightarrowfill}}
\limits^{\textstyle#1}_{\textstyle#2}}}

\def\lfl#1#2{\smash{\mathop{\hbox to 12mm{\leftarrowfill}}
\limits^{\textstyle#1}_{\textstyle#2}}}

\def\dfl#1#2{\llap{$\textstyle#1$}\left\downarrow
\vbox to 6mm{}\right.\rlap{$\textstyle#2$}}

\def\rfll#1#2#3{\smash{\mathop{\hbox to#3mm{\rightarrowfill}}
\limits^{\textstyle#1}_{\textstyle#2}}}

\def\lfll#1#2#3{\smash{\mathop{\hbox to#3mm{\leftarrowfill}}
\limits^{\textstyle#1}_{\textstyle#2}}}

\def\diagram#1{\def\normalbaselines{\baselineskip=0pt
\lineskip=10pt\lineskiplimit=1pt}\matrix{#1}}

%%%%%%%%%%%%%%%%%%%%%%%%%%%%%%%%%%%%%%%%%%%%%%%%%%%%%%%%%%%%%%%%%%%%%%%%%%%
%% THIS IS forum.tex A MACRO FILE FOR TeX USERS OF Semigroup Forum
%% VERSION OF 1992-02-29.  
%% Dimensions revised, firstpage improved, and macros perfected.
%% AUTHOR: SAVE THIS FILE AS forum.tex IN A DIRECTORY WHICH TeX SEARCHES FOR 
%% \input FILES.
%%%%%%%%%%%%%%%%%%%%%%%%%%%%%%%%%%%%%%%%%%%%%%%%%%%%%%%%%%%%%%%%%%%%%%%%%%%
\magnification=1200
\tolerance=300
\pretolerance=200
\hfuzz=1pt
\vfuzz=1pt
\parindent=35pt
\mathsurround=1pt
\parskip=1pt plus .25pt minus .25pt
\baselineskip=11pt
\normallineskiplimit=1pt
\abovedisplayskip=10pt plus 2.5pt minus 7.5pt
\abovedisplayshortskip=0pt plus 2.5pt
\newdimen\jot
\jot=2.5pt

\hoffset=0.3in          
\voffset=-1\baselineskip 
\hsize=5.8 true in
\vsize=53.25\baselineskip  

\font\eightrm=cmr8

\font\bfone=cmbx10 scaled 1200
\font\smc=cmcsc10

%%% If you have the AMS-TeX fonts msam and msbm, then uncomment the
%%% following lines
%\font\tenmsc=msbm10 at 12pt
%\font\sevenmsc=msbm8 
%\font\fivemsc=msbm7
%\font\tenmsb=msbm10 
%\font\sevenmsb=msbm7
%\font\fivemsb=msbm5
%\font\tenmsa=msam10 
%\font\sevenmsa=msam7
%\font\fivemsa=msam5
%\newfam\msbfam
%\newfam\msafam
%\newfam\mscfam
%\textfont\mscfam=\tenmsc  \scriptfont\mscfam=\sevenmsc
%  \scriptscriptfont\mscfam=\fivemsc
%\textfont\msbfam=\tenmsb  \scriptfont\msbfam=\sevenmsb
%  \scriptscriptfont\msbfam=\fivemsb
%\textfont\msafam=\tenmsa  \scriptfont\msafam=\sevenmsa
%  \scriptscriptfont\msafam=\fivemsa
%
%\def\hexnumber@#1{\ifcase#1 0\or1\or2\or3\or4\or5\or6\or7\or8\or9\or
%        A\or B\or C\or D\or E\or F\fi }
%
%\edef\msb@{\hexnumber@\msbfam}
%\edef\msa@{\hexnumber@\msafam}
%\edef\msc@{\hexnumber@\mscfam}
%%% We do not use AMS-TeX -----hence
%\def\Bbb{\ifmmode\let\next\Bbb@\else
% \def\next{\errmessage{Use \string\Bbb\space only in math mode}}\fi\next}
%\def\Bbb@#1{{\Bbb@@{#1}}}
%\def\Bbb@@#1{\fam\mscfam#1}
\mathchardef\emptyset="001F

\def\nin{\noindent}

\nopagenumbers
\def\rightheadline{\hfil\smc\lastname\hfil\tenbf\folio}
\def\leftheadline{\tenbf\folio\hfil\smc\lastname\hfil}
\def\firstheadline{\vbox{\baselineskip=8pt
    \parindent 0pt \obeylines \eightrm
    Semigroup Forum  \hfill{\date}
    \par\nin Springer-Verlag New York Inc.}}
\headline={\ifnum\pageno=1\firstheadline
           \else\ifodd\pageno\rightheadline
           \else\leftheadline\fi\fi}
    %IN ORDER TO ELIMINATE THE SOURCE DATA ON THE FIRST TWO LINES OF
    %THE FIRST PAGE FROM A PRINTOUT, TYPE INTO THE FINAL DOCUMENT
    %THE LINE \def\firstheadline{\vbox{\vskip2baselineskip}}
\def\datum{\ifcase\month\or January\or February\or March\or April
    \or May\or June\or July\or August\or September\or October
    \or November\or December\fi\space\number\day, \number\year}
\def\date{{\eightrm Version of~}\datum}
    %IN ORDER TO ELIMINATE THE VERSION DATE
    %FROM THE FINAL PRINTOUT, TYPE INTO THE FINAL DOCUMENT
    %THE LINE  \def\date{}

\def\title{AUTHOR: INSERT TITLE!!}
\def\author{AUTHOR: INSERT YOUR NAME(S)}
\def\lastname{AUTHOR: INSERT YOUR LASTNAME(S) (WITHOUT INTIALS)}
\def\editor{AUTHOR: INSERT EDITOR'S NAME}
\def\thanks#1{\footnote*{\eightrm#1}}

\def\title{\centerline{\bfone\titleone}\centerline{\bfone\titletwo}}
\def\titleone{}
\def\titletwo{}    
\def\firstpage{\vbox{\vskip1.5truecm}
       \title 
       \vskip.75truecm\centerline{\bf\author}
       \medskip
       \centerline{\eightrm Communicated by \editor}
       \vskip.75truecm\rm}

\def\sectionheadline #1{\vskip-\lastskip\bigbreak
                      \centerline{\bf #1}\nobreak\medskip\nobreak}

        \newtoks\literat
\def\[#1 #2\par{\literat={#2\unskip.}%
        \hbox{\vtop{\hsize=.1\hsize\nin [#1]\hfill}
        \vtop{\hsize=.9\hsize\nin\the\literat}}\par
        \vskip.3\baselineskip}

\def\references{\sectionheadline{\bf References}\frenchspacing
                \entries\par}

\def\address{Author: {\tt$\backslash$def$\backslash$address$\{$??$\}$}}
\def\addresstwo{}  
\def\lastpage{\par\vbox{\vskip.8truecm\nin
\line{\vtop{\hsize=.45\hsize{\eightrm\parindent=0pt\baselineskip=8pt
            \nin\address}}\qquad
      \vtop{\hsize=.45\hsize\nin{\eightrm\parindent=0pt
            \baselineskip=8pt\addresstwo}}\hfill}
}}

\def\qed{{\unskip\nobreak\hfil\penalty50\hskip .001pt \hbox{}\nobreak\hfil
          \vrule height 1.2ex width 1.1ex depth -.1ex
          \parfillskip=0pt\finalhyphendemerits=0\medbreak}\rm}
  %THIS IS THE END-OF-PROOF SIGN.
  %NOT TO BE USED IN DISPLAY MODE! IF YOU WANT TO CONCLUDE A PROOF
  %AT THE END OF A LINE IN DISPLAY MODE USE
\def\qeddis{\eqno{\vrule height 1.2ex width 1.1ex depth -.1ex} $$
                   \medbreak\rm}
  %BUT OMIT $$---THE MACRO WILL WRITE THAT!

\def\Proposition #1. {\bigbreak\vskip-\parskip\noindent{\bf Proposition #1.}
    \quad\it}
\def\Propositions. {\bigbreak\vskip-\parskip\noindent{\bf Proposition.}
    \quad\it}    
\def\Theorem #1. {\bigbreak\vskip-\parskip\noindent{\bf  Theorem #1.}
    \quad\it}
\def\Corollary #1. {\bigbreak\vskip-\parskip\nin{\bf Corollary #1.}
    \quad\it}
\def\Lemma #1. {\bigbreak\vskip-\parskip\noindent{\bf  Lemma #1.}\quad\it}
\def\Lemmad. {\bigbreak\vskip-\parskip\noindent{\bf  Lemma.}\quad\it}
\def\Definition #1. {\rm\bigbreak\vskip-\parskip\noindent{\bf Definition #1.}
    \quad}
\def\Remark #1. {\rm\bigbreak\vskip-\parskip\noindent{\bf Remark #1.}\quad}
\def\Remarkd. {\rm\bigbreak\vskip-\parskip\noindent{\bf Remark.}\quad}
\def\Exercise #1. {\rm\bigbreak\vskip-\parskip\noindent{\bf Exercise #1.}
    \quad}
\def\Example #1. {\rm\bigbreak\vskip-\parskip\noindent{\bf Example #1.}\quad}
\def\Examples #1. {\rm\bigbreak\vskip-\parskip\noindent{\bf Examples #1.}\quad}
\def\Proof#1.{\rm\par\ifdim\lastskip<\bigskipamount\removelastskip\fi
    \smallskip\noindent{\bf Proof.}\quad}
%end of forum.tex
%%%%%%%%%%%%%%%%%%%%%%%%%%%%%%%%%%%%%%%%%%%%%%%%%%%%%%%%%%%%%%%%%%%%%%%%%%%%%

%% THE USER MAY WISH TO REDEFINE \hoffset AND/OR \voffset TO
%% SUIT LOCAL NEEDS

\hoffset=0.3in
\voffset=-1\baselineskip

%% IF YOU USE US-LETTERFORMAT FOR YOUR HARDCOPY PRINTOUTS YOU MAY HAVE TO
%% CHANGE

\vsize=55\baselineskip

%% TO

%\vsize=54\baselineskip

%% INSERT THE FIRST LINE OF YOUR TITLE BETWEEN "{" AND "}"

\def\titleone{Embedding simple commutative monoids}

%% IF THERE IS A SECOND LINE TO YOUR TITLE INSERT IT BETWEEN "{" AND "}"

\def\titletwo{into simple refinement monoids}

%% INSERT YOUR NAME(S) BETWEEN "{" AND "}"

\def\author{Friedrich Wehrung}

%% INSERT YOUR LAST NAME(S) BETWEEN "{" AND "}"; IN CASE OF MORE THAN ONE
%% AUTHOR, SEPARATE THE LASTNAMES BY ``and"

\def\lastname{Wehrung}

%% INSERT THE NAME OF THE EDITOR TO WHOM YOU SUBMIT THE PAPER
%% BETWEEN "{" AND "}".

\def\editor{Klaus Keimel}

%% INSERT YOUR REFERENCES BETWEEN "{" AND "}", AS IN THE SAMPLE BELOW.
%% NOTE THAT INDIVIDUAL ENTRIES DO NOT TERMINATE WITH A PERIOD, BUT
%% INSTEAD WITH A BLANK LINE. THIS APPLIES TO THE LAST ENTRY AS WELL.

\def\entries
{\[1 Bergman, G. M., {\it unpublished correspondence} (December 1990)

\[2 Clifford, A. H., and G. B. Preston, ``The algebraic theory of
semigroups,'' Mathematical Surveys {\bf7}, American Mathematical Society,
Providence,  R.I., vol. {\bf1}, 1961 and vol. {\bf2}, 1967

\[3 Dobbertin, H., {\it Refinement monoids, Vaught Monoids, and
Boolean Algebras}, Mathematische Annalen {\bf 265} (1983),
475--487

\[4 Goodearl, K. R., ``Von Neumann Regular Rings'', Pitman, London, 1979

\[5 Goodearl, K. R., ``Partially ordered abelian groups with the
interpolation pro\-per\-ty'', Mathematical surveys and monographs,
number {\bf 20}, American Mathematical Society, 1986

\[6 Goodearl, K. R., {\it Von Neumann regular rings and direct sum
decomposition problems}, Abelian Groups and Modules, Padova
1994 (A. Facchini and C. Menini, eds.), Dordrecht (1995) Kluwer,
249--255

\[7 Grillet, P. A., {\it Interpolation properties and tensor products
of semigroups}, Semigroup Forum {\bf 1} (1970), pp. 162-168

\[8 Howie, J. M., ``An introduction to semigroup theory'', Academic
Press: London, New York, San Francisco, 1976

\[9 Moncasi, J., {\it A regular ring whose $K_0$ is
not a Riesz group}, Communications in Algebra {\bf 13}
(1985), 125--133

\[10 Pardo, E., {\it On the representation of simple Riesz groups},
to appear in Communications in Algebra

\[11 Wehrung, F., {\it Injective positively ordered monoids I},
Journal of Pure and Applied Algebra {\bf 83} (1992), 43--82

\[12 Wehrung, F., {\it Restricted injectivity, transfer property
and decompositions of separative positively ordered monoids},
Communications in Algebra {\bf 22} (5) (1994), 1747--1781

\[13 Wehrung, F., {\it Monoids of intervals of ordered abelian groups},
Journal of Algebra {\bf 182} (1996), pp. 287--328

\[14 Wehrung, F., {\it Tensor products of structures with
interpolation},
Pacific Journal of Mathematics {\bf 176} (1) (1996), 267--285

}

% \def\entries
% {\[1 Chopin, F., and G. Sand, {\it Harmonic analysis on a Steinway
% manifold,}
% Proc. Acad. Nat. Polonaise {\bf 43}(1832), 143--169
%
% \[2 Hilgert, J., and K. H. Hofmann, {\it Invariant cones in Lie
% algebras,}
% Semigroup Forum {\bf 37}(1988), 241--252
%
% \[3 Knuth, Donald, E., ``The \TeX{book},'' Addison-Wesley,
% Reading, Mass., 1984
%
% }

%% INSERT YOUR ADDRESS BETWEEN "{" AND "}"; LEAVE A BLANK LINE BETWEEN
%% EACH LINE OF THE ADDRESS BUT DO NOT START WITH AN EMPTY LINE.
%% ADDRESSES ARE WRITTEN WITHOUT REPEATING THE AUTHOR'S NAME.
%% ALSO, TERMINATE THE ADDRESS WITH A BLANK LINE

\def\address{
Universit\'e de Caen

D\'epartement de Math\'ematiques

14032 Caen Cedex

FRANCE

e-mail: gremlin@math.unicaen.fr

}

%% PLEASE FILL OUT IF THERE IS A SECOND AUTHOR

\def\addresstwo{

}

%% BEGIN WRITING YOUR TEXT AFTER "\firstpage"

\firstpage

\footnote{}{\sevenrm 1991 {\sevensl Mathematics subject
classification}:
\sevenrm 06F20, 20M14.\par
{\sevensl Key words and phrases}:
\sevenrm refinement monoids; simple \cm s;
order-embeddings; unitary embeddings; di\-vi\-si\-bi\-li\-ty.}

{\sevenbf Abstract.} {\sevenrm Say that a {\sevenit cone} is a
\cm\ that is in addition {\sevenit conical}, {\sevenit i.e.},
satisfies $\scriptstyle x+y=0\ \Rightarrow\ x=y=0$.
We show that cones (resp.
simple cones) of many kinds order-embed or even embed
unitarily into refinement cones (resp. simple refinement
cones) of the same kind, satisfying in addition various
divisibility conditions. We do this in particular for all cones,
or for all separative cones, or for all cancellative cones
(positive cones of \poag s). We also settle both the
torsion-free case and the unperforated case. Most of our
results extend to arbitrary \cm s.}

\bigskip

\sectionheadline{\S0. Introduction.}

In [{\bf\goodII}, Question 30], K. R. Goodearl asks whether
every simple \poag\ can be embedded into a simple Riesz group,
and similarly for simple torsion-free (resp. unperforated)
\poag s (the condition ``simple" has to be added among the
hypotheses --- indeed, every subgroup of a simple \poag\ is
simple). In this paper we solve positively this
question, and many others of the same ilk concerning {\it
cones}, {\it i.e.}, conical \cm s
(not necessarily cancellative). A prototype can be found in
P. A. Grillet's [{\bf\gril}, Theorem 1],
where the author outlines an embedding procedure of every \cm\
into a refinement monoid; the same result, with a more
detailed (and probably independent) proof can be found in H.
Dobbertin's [{\bf\dobb}, Theorem 5.1]; {\it see} also
[{\bf\wehrII}] for more general results where the
{\it preordering} is also taken into consideration.
The Grillet-Dobbertin embedding result was improved in
[{\bf\berg}] by G. M. Bergman, who essentially proved there a
version of our Lemma 1.6, which implies in particular that
every \cm\ can be order-embedded into a divisible refinement
monoid (in his notes, Bergman uses the notion of ``saturated
embedding", which is in fact equivalent to the classical
notion of unitary embedding [{\bf\howi}] and thus equivalent
to the definition of unitary embedding adopted in this paper
minus cofinality). Using this, Bergman proved also that every
simple conical \cm\ with an infinity ({\it i.e.}, an element
$\infty$ satisfying $(\forall x)(x+\infty=\infty)$) embeds
into a simple conical refinement monoid with an infinity.
Then, using submonoids of products of these with ${\bb R}^+$, he
produced examples of non-cancellative directly finite simple
refinement monoids, thus solving a problem of Goodearl and
O'Meara. Let us emphasize that Bergman's notes [{\bf\berg}]
contain several other interesting results, and, although our
work is related to his, it has been done independently.
\smallskip

Throughout this paper, we will systematically consider
additional structural information about both the monoids
(separativity, simplicity...) and the embeddings
(order-embedding, unitarity, strong unitarity). Thus most of
our results will be of the type ``{\it every cone of such and
such kind admits an embedding of such and such kind into a cone
satisfying such and such closure property}", and ``{\it every
cone with such and such closure property is (among other
properties) a refinement monoid}". Let us be more
specific:\smallskip

\item{---} ``{\it Every cone of such and such kind...}".
Some of the classes of cones that we will consider will be the
class of all cones, or the class of all {\it cancellative}
(resp. {\it separative}, {\it stably finite},
{\it unperforated}...) cones.\smallskip

\item{---} ``{\it ...an embedding of such and such kind...}".
Three sorts of embeddings will be considered: {\it
order-embeddings} ({\it i.e.}, monoid homomorphisms that are in
addition embeddings for the algebraic preorderings on each
side), and the stronger very important notion of {\it unitary
embeddings} or {\it strong unitary embeddings}
(Definition 1.2) which will allow us to transfer
the embedding results to, {\it e.g.}, {\it separative} cones
(it can also been shown, in the language of universal algebra,
that any unitary embedding has the ``congruence extension
property").\smallskip

\item{---} ``{\it ... such and such closure property}". These
closure properties will always be {\it maximality} properties
(similar to algebraic closure)
for existence of solutions of finite equation systems
(Definition 1.4) in extensions of the monoid under
consideration. Their main advantage will be that their models
satisfy several properties at once, as refinement, divisibility
or quasi-divisibility, {\it etc.}.\medskip

We now summarize briefly each section.\smallskip

\item{---} In section 1, we will prove embedding theorems for
cones, without considering simplicity yet. Proposition
1.5 states that cones from a given class embed (sometimes
unitarily) into ``closed" cones; its proof is in fact a
general argument, reminiscent (with a similar proof) of the
model-theoretical fact that every model of a
$\forall\exists$ theory embeds into an existentially closed
model of that theory (but the latter does not take into
account unitary embeddings, nor simplicity). Theorems 1.8 and
1.9 will yield the fact that among other properties, ``closed"
cones are refinement monoids. This works as well for all cones
as, {\it e.g.}, for cancellative, separative, stably finite,
{\it etc.}, cones. Theorem 1.14 extends this to
torsion-freeness or unperforation, but then the following
phenomenon happens: {\it either the embedding has very strong
properties and we only keep refinement, or the embedding is
nothing more than an order-embedding but we have both
refinement and divisibility}. In a few words,\smallskip

\centerline{\bfi What is gained for the embedding is lost}
\centerline{\bfi for the final object, and vice-versa.}
\smallskip

\item{\phantom{---}}This phenomenon will be
present in all our results, and explains their great variety,
which forced us to restrict ourselves to a sample. {\it In
particular, most of our results extend to the class of all
\cm s, with very similar proofs}, but we chose not to give
details about this in order to maintain this paper down to a
reasonable size; moreover, the class of \cm s which we had in
mind often arise from K-theoretical situations, and these are
always conical. Proposition
1.16 shows that ``closed" cones satisfy the rather mysterious
axiom (WSD), which in turn, together with the finite refinement
property, guarantees other forms of refinement (as, {\it e.g.},
the finite interpolation property, or the ``hereditary
refinement" property (HREF) both considered in
[{\bf\wehrIII}]).\smallskip

\item{---} Section 2 extends the results of section 1 to
{\it simple} cones. Proposition 2.5 is the ``simple analogue"
of Proposition 1.5 (existence of order-embeddings or of
unitary embeddings into ``simply closed" cones), while the
``simple analogue" of Theorem 1.9 is Theorem 2.6 (``simply
closed" cones are, among other things, refinement
cones). It is also in this section that are
obtained answers to Goodearl's question for either arbitrary
or unperforated \poag s (Corollaries 2.8 and 2.9), {\it but
not for torsion-free \poag s yet}. Finally, Proposition 2.12
shows that ``simply closed" cones satisfy (WSD).\smallskip

\item{---} Section 3 is devoted to the remaining {\it
torsion-free} case, which requires a special treatment because
of problems of divisibility --- thus the trick of using
instead the property of {\it quasi-divisibility} (Definition
3.2). We will in fact more generally
consider torsion-freeness with respect to a given non-trivial
multiplicative subsemigroup ${\bb P}$ of ${\bb N}$. Then most
results of previous sections extend to this context ---
yielding essentially Theorems 3.5 and 3.9. In
particular, the remaining problem about embedding
simple torsion-free \poag s into simple torsion-free Riesz
groups is solved in Corollary 3.12.\medskip

We shall put
${\bb N}={\bb Z}^+\setminus\{0\}$. If $a_i$, $b_j$, $c_{ij}$ ($i,j<2$)
are elements of some \cm\ $M$, then we will say that the
following array

$$\vbox{\offinterlineskip
\halign{\tv\cc{$#$}&\tv\cc{$#$}&\tv\cc{$#$}&\tv#\cr
\omit&\omit\hrulefill&\omit\hrulefill&\omit\cr
\omit\tvi&b_0&b_1&\cr
\noalign{\hrule}
a_0&c_{00}&c_{01}&\cr
\noalign{\hrule}
a_1&c_{10}&c_{11}&\cr
\noalign{\hrule}
}}$$

\noindent is a ($2\times2$) {\it refinement matrix} when for
all $i<2$, we have $a_i=c_{i0}+c_{i1}$ and $b_i=c_{0i}+c_{1i}$.
If such a matrix always exists provided that $a_0+a_1=b_0+b_1$,
then we will say that $M$ is a {\it refinement monoid}.
If $M$ is a \cm, the {\it algebraic} preordering
$\leq_{\rm alg}$ of $M$ is defined by
$x\leq_{\rm alg}y\Leftrightarrow(\exists z)(x+z=y)$
(this terminology is borrowed from the ${\rm C}^*$-algebraists
and we find it more inspired than the term ``minimal" used in
[{\bf\wehr}, {\bf\wehrII}], thus we will use it throughout
this text). We will often identify a \cm\ $M$ with its adjunct
$(M,\leq_{\rm alg})$, emphasized by the denomination
``algebraic \cm". The notion of homomorphism is not affected by
this adjunction, contrarily to the notion of embedding; thus if
$M$ and $N$ are \cm s, an {\it order-embedding} $f$ from $M$
into $N$ will be an injective homomorphism of monoids which
satisfies in addition $x\leq y\Leftrightarrow f(x)\leq f(y)$.
If $M$ is a submonoid of a \cm\ $N$, we will say that $N$ is
an {\it extension} of $M$ when the inclusion map from $M$ into
$N$ is an order-embedding.

An {\it order-unit} of a
\cm\ $M$ is an element $u$ of $M$ such that
$(\forall x\in M)(\exists n\in{\bb N})(x\leq_{\rm alg}nu)$; $M$ is
{\it simple} when every $u\in M$ such that $u>0$ ({\it i.e.},
$u\not\leq_{\rm alg}0$) is an order-unit of $M$, {\it conical}
when it satisfies $(\forall x,y)(x+y=0\Rightarrow x=y=0)$,
{\it cancellative} when it satisfies
$(\forall x,y,z)(x+z=y+z\Rightarrow x=y)$,
{\it separative} [{\bf\clpr}, vol. 1] when it satisfies 
$(\forall x,y)(2x=x+y=2y\Rightarrow x=y)$. Note that the
latter notion is (strictly) weaker than the notion of
separativity used in [{\bf\wehrII}], which was designed there
for (positively) {\it preordered} monoids.\smallskip

In addition, if $M$ and $N$ are two \cm s and $f:M\to N$ is a
homomorphism of monoids, we will say that $f$ is {\it conical}
when $f^{-1}\{0_N\}=\{0_M\}$ (note that we write
$f^{-1}\{0\}$ instead of ${\rm Ker}(f)$, since in
monoid-theoretical contexts, the latter denotes in
general rather a {\it congruence} than a {\it submonoid}).
Thus if $f$ is one-to-one, then it is conical, but the
converse is false.\smallskip

Furthermore, we will use the canonical transpositions of
definitions used for regular rings
[{\bf\good}, {\bf\goodIII}] to the monoid world: thus,
$M$ will be said to be {\it stably finite} when it satisfies
the axiom $(\forall x,y)(x+y=y\Rightarrow x=0)$. We finally
refer to [{\bf\goodII}] for the terminology about \poag s.
\bigskip\goodbreak

\sectionheadline{\S1. Order-embedding cones.}

In this chapter we will embed algebraic \cm s,
without considering simplicity yet. In fact, most of our
results will concern {\it conical} \cm s:\medskip

\Definition1.1. A {\it cone} is a conical
\cm; a {\it refinement cone} is a cone satisfying
the finite refinement property.\medskip

We will need some algebraic theory of \cm s, in particular the
study of {\it amalgamated sums}. For this purpose, 
we will slightly deviate
from the original notion of {\it unitary embedding}
[{\bf\howi}, page 232] used in the context of 
(non-commutative) semigroup amalgams, with the adjunction of
a {\it cofinality} condition:\medskip

\Definition1.2. Let $f:A\to B$ be a
homomorphism of \cm s. Then $f$ is {\it unitary} when it
is {\it one-to-one}, has {\it cofinal} image (for the
algebraic preordering of $B$) and it satisfies
$$(\forall a_0,a_1\in A)
(\forall b\in B)\bigl(f(a_0)+b=f(a_1)\Rightarrow
b\in f[A]\bigr).$$
We will say that $f$ is {\it strongly unitary}
when it is unitary and satisfies, for all $m\in{\bb N}$, the
condition
$$(\forall b\in B)
\bigl(mb\in f[A]\Rightarrow b\in f[A]\bigr).$$
In the case where $A$ is a submonoid of $B$
and $f$ is the inclusion map, we will say that
$B$ is a {\it unitary extension} (resp. a {\it strong unitary
extension}) of $A$.\medskip

Important examples of unitary and strong unitary embeddings are
the following: let $A$ be a cofinal subgroup of a \poag\ $B$.
Then $B^+$ is a unitary extension of $A^+$.
If in addition $B/A$ is {\it torsion-free}, then $B^+$ is a
strong unitary extension of $A^+$; the converse holds provided
for example that $A$ is {\it directed}.
\smallskip

The proof of the following lemma is straightforward:\medskip

\Lemma1.3.
\item{{\rm(a)}} Every unitary homomorphism of \cm s is an
order-em\-bed\-ding.
\item{{\rm(b)}} If $f:A\to B$ and $g:B\to C$ are unitary
(resp. strongly unitary), then $g\circ f$ is unitary
(resp. strongly unitary).
\item{{\rm(c)}} If $B$ is a unitary (resp. strong unitary)
extension of $A$ and $C$
is a monoid such that $A\subseteq C\subseteq B$, then $C$ is
a unitary (resp. strong unitary) extension of $A$.
\item{{\rm(d)}} If $A$ is a \cm, then any directed union of
unitary (resp. strong unitary) extensions of $A$ is a unitary
(resp. strong unitary) extension of $A$.\qed\medskip

If $M$ is a \cm, say that an {\it equation system} with
parameters from $M$ is a finite set of ``equations" 
each of them of the form
$$\sum_{i<k}p_i{\bf x}_i+a=\sum_{i<k}q_i{\bf x}_i+b$$
where $k$ and $p_i$, $q_i$ ($i<k$) are non-negative integers,
$a$, $b$ are elements of $M$ (``parameters") and
${\bf x}_i$ ($i<k$) are symbols of variable (``unknowns"). The
notion of {\it solution} of a given equation system should
then be obvious. We will need in this paper three
definitions, related to equation systems and unitarity:\medskip

\Definition1.4. Let {\bf C}
be a class of \cm s and let $M$ belong to {\bf C}. Then $M$ is
said to be {\it{\bf C}-closed} (resp.
{\it{\bf C}-unitarily closed},
{\it{\bf C}-strongly unitarily closed}) when for every
equation system $\Sigma$ with parameters from $M$, if $\Sigma$
admits a solution in some extension (resp. some
unitary extension, resp. some strong unitary
extension) of $M$ belonging to {\bf C}, then it admits a
solution in $M$.\medskip

Note the obvious {\it reverse} implications
$$\hbox{(strongly unitarily
closed)}\Leftarrow\hbox{(unitarily closed)}\Leftarrow
\hbox{(closed)}.$$
It would perhaps be preferable, in a more general context, to
use the terminology ``unitarily {\it algebraically} closed",
{\it etc.}, but we choose not to do so for simplicity
sake.\smallskip

If {\bf C} is an arbitrary class of \cm s, to be ``simply
{\bf C}-closed" will be by definition the same as to be
${\bf C}'$-closed where ${\bf C}'$ is the class of all simple
elements of {\bf C}. Similar conventions will apply to ``simply
{\bf C}-unitarily closed", {\it etc.}. A general existence
theorem is then the following:\medskip

\Proposition1.5. Let {\bf C} be a class of
\cm s closed under unions of chains. Then every element of
{\bf C} admits an order-embedding (resp. a unitary embedding, a
strong unitary embedding) into a {\bf C}-closed (resp. a
{\bf C}-unitarily closed, {\bf C}-strongly unitarily closed)
element of {\bf C}.
\Proof. Obtained by aping the classical proof of
embedding any model of a theory with class of models closed
under unions of chains into an existentially closed model: let
$M$ be an element of {\bf C} and let
$(\varphi_\xi)_{\xi<\kappa}$ be an enumeration of all equation
systems with parameters from $M$ (where $\kappa$ is an
infinite cardinal). Using Lemma 1.3 (d) and the fact that
{\bf C} is closed under unions of chains, it is not
difficult to construct a transfinite chain
$(M_\xi)_{\xi<\kappa}$ of elements of {\bf C} satisfying the
following properties:\smallskip

\item{{\rm(i)}} If $\xi<\eta<\kappa$, then $M_\eta$ is an extension
(resp. a unitary extension, a strong unitary extension) of
$M_\xi$.

\item{{\rm(ii)}} For all $\xi<\kappa$, if $\varphi_\xi$ admits a
solution in some extension (resp. unitary extension,
strong unitary extension) of $M_\xi$ in {\bf C},
then it admits a solution in $M_{\xi+1}$.\smallskip

Then put $M'=\bigcup_{\xi<\kappa}M_\xi$. Then $M'\in{\bf C}$
since {\bf C} is closed under unions of chains, and $M'$ is an
extension (resp. a unitary extension, a
strong unitary extension) of $M$. Then define
$M^{(n)}$ ($n\in{\bb Z}^+$) by $M^{(0)}=M$, and
$M^{(n+1)}=(M^{(n)})'$ for all $n\in{\bb Z}^+$. Then
$\tilde M=\bigcup_{n\in{\bb Z}^+}M^{(n)}$ satisfies the required
conditions.\qed\medskip

Thus the largest part of this paper will be to investigate the
structure of all {\bf C}-closed (resp.
{\bf C}-unitarily closed, {\bf C}-strongly unitarily closed)
\cm s for various classes {\bf C}. Thus it is natural to start
with {\bf C} being the class of all \cm s.\medskip

\Lemma1.6. Let $e:A\to B$ and $f:A\to C$
be homomorphisms of \cm s. Let $D=B\amalg_{e,f}C$ be the
amalgamated sum of $B$ and $C$ along $e$ and $f$ in the
category of \cm s, so that we have a commutative diagram
\goodbreak
$$\diagram{
A&\rfl{f}{}&C\cr
\dfl{e}{}&&\dfl{\bar e}{}\cr
B&\rfl{\bar f}{}&D\cr
}$$
Then the following holds:\smallskip
\item{{\rm(a)}} If both $B$ and $C$ are
conical and if both $e$ and $f$ are conical maps, then
$D$ is conical and both $\bar e$ and $\bar f$ are conical maps
(thus $D$ is also the amalgamated sum of $B$ and
$C$ along $e$ and $f$ in the category of cones with conical
homomorphisms).
\item{{\rm(b)}} If $e$ is unitary, then $\bar e$ is also
unitary.
\item{{\rm(c)}} If $e$ is strongly unitary, then
$\bar e$ is also strongly unitary.
\smallskip
{\rm In particular, ``unitarity and strong unitarity
are transferable".}
\Proof. Let $\to$ be the binary relation on
$B\times C$ consisting exactly on the pairs
$$(b+e(a),c)\to(b,f(a)+c)\qquad({\rm all}\ a\in A,\ b\in B\ 
{\rm and}\ c\in C).$$
Then $\to$ is compatible with the addition ({\it i.e.},
$\xi\to\eta$ implies $\xi+\zeta\to\eta+\zeta$), thus so is the
equivalence relation $\equiv$ on $B\times C$ generated by
$\to$. It follows that $D=B\times C/\!\!\equiv$ (up to a
natural isomorphism). For all $(b,c)\in B\times C$, denote by
$[b,c]$ its natural image in $D$. Then the natural maps
$\bar e$ and $\bar f$ are defined by $\bar e(c)=[0,c]$ and
$\bar f(b)=[b,0]$.\smallskip

Let us first settle (a). To prove that $D$ is conical
and that $\bar e$ and $\bar f$ are conical, it
suffices, since both $B$ and $C$ are conical, to prove that
for all $(b,c)\in B\times C$, $[b,c]=[0,0]$ implies that
$(b,c)=(0,0)$ (the converse being trivial); thus it suffices
in fact to prove the conclusion for $(b,c)\to(0,0)$ and
$(0,0)\to(b,c)$. In the first case, there exists $a\in A$ such
that $b=0+e(a)$ and $0=f(a)+c$; by assumption on $C$ and
$f$, it follows that $(b,c)=(0,0)$. In the second case, there
exists $a\in A$ such that $0=b+e(a)$ and $c=f(a)+0$; by
assumption on $B$ and $e$, it follows again that
$(b,c)=(0,0)$. The conclusion of (a) follows.\smallskip

Suppose from now on that $e$ is unitary. Without loss of
generality, $A\subseteq B$ and $e$ is the inclusion map from
$A$ into $B$, supposed to be unitary.

Since $A$ is cofinal in $B$, for all 
$(b,c)\in B\times C$, there exists $a\in A$ such that
$b\leq a$, thus $[b,c]\leq[a,c]=[0,f(a)+c]=\bar e(f(a)+c)$, so
that the image of $\bar e$ is cofinal in $D$.
\medskip

\noindent{\bf Claim.} {\sl For all $a\in A$, $b\in B$ and
$c,c'\in C$, if $(a,c)\equiv(b,c')$, then $b\in A$ and
$f(a)+c=f(b)+c'$.}\smallskip

\noindent{\bf Proof of Claim.} It clearly suffices to show
that the conclusion holds provided that either
$(a,c)\to(b,c')$ or $(b,c')\to(a,c)$. In the first case,
there exists $a'\in A$ such that $a=b+a'$ and $c'=f(a')+c$;
since $e$ is unitary, $b\in A$ and
$f(a)+c=f(b)+f(a')+c=f(b)+c'$. In the second case, there
exists $a'\in A$ such that $b=a+a'$ and $c=f(a')+c'$; thus
$b\in A$ and $f(a)+c=f(a)+f(a')+c'=f(b)+c'$.\qed\medskip

From the Claim above, it follows immediately that $\bar e$ is
one-to-one. Furthermore, if $c_0,c_1\in C$ and
$(b,c)\in B\times C$ such that 
$\bar e(c_0)+[b,c]=\bar e(c_1)$, {\it i.e.},
$(b,c_0+c)\equiv(0,c_1)$, then, by the Claim above, $b\in A$,
whence $[b,c]=\bar e(f(b)+c)\in\bar e[C]$, which proves that
$\bar e$ is unitary; thus (b).\smallskip

Finally suppose that $e$ is strongly unitary. Let
$m\in{\bb N}$ and let $[b,c]\in D$ such that
$m\cdot[b,c]\in\bar e[C]$. This means that there exists
$c'\in C$ such that $[mb,mc]=[0,c']$. By the Claim above,
$mb\in A$; since $B$ is a strong unitary extension of $A$, it
follows that $b\in A$. Therefore,
$[b,c]=[0,f(b)+c]\in\bar e[C]$. Thus $\bar e$ is strongly
unitary, and (c) follows.
\qed\medskip

Now let $R=\{(x_0,x_1,x_2,x_3)\in{\bb Z}^4:\ x_0+x_3=x_1+x_2\}$,
equipped with the ordering induced by the natural ordering of
${\bb Z}^4$.\medskip

\Lemma1.7.
\item{{\rm(a)}} {\sl $R^+$ is the submonoid of
$({\bb Z}^+)^4$ generated by $\alpha_0=(1,1,0,0)$,
$\alpha_1=(0,0,1,1)$, $\beta_0=(1,0,1,0)$ and
$\beta_1=(0,1,0,1)$.}
\item{{\rm(b)}} {\sl For every \cm\ $M$ and all elements $a_i$,
$b_i$ ($i<2$) of $M$ such that $a_0+a_1=b_0+b_1$, there exists
a unique homomorphism $f:R^+\to M$ such that for all $i<2$,
$f(\alpha_i)=a_i$ and $f(\beta_i)=b_i$. If in addition $M$ is
conical and $a_i,b_i\ne0$ (all $i<2$), then
$f$ is conical.}
\Proof. Most of it is well-known --- {\it see}
for example [{\bf\wehrII}, proof of Corollary 2.7]. The last
part of (b) is trivial.\qed\medskip

We now come to the main theorem of this section (whose full
meaning is highlighted by Proposition 1.5):\medskip

\Theorem1.8.
\item{{\rm(a)}} Every strongly unitarily
closed cone is a refinement cone.\smallskip
\item{{\rm(b)}} Let $M$ be a unitarily closed cone. Then
$M$ is ``{\rm normally divisible}", {\rm i.e.},
for every $p\in{\bb N}$, every finite subset $X$ of
$M$ and every element $a$ of $M$, there exists $u\in M$ such
that $pu=a$ and for all $x,y\in X$, $x+a=y+a$ implies
$x+u=y+u$.\smallskip
{\rm In the statement of (a), a ``strongly unitarily closed cone"
is by definition a cone that is strongly unitarily closed
relatively to the class of all cones, and similarly for
(b).}
\Proof. We start with (a). Let $a_i$, $b_i$
($i<2$) be elements of $M$ such that $a_0+a_1=b_0+b_1$.
Consider the following equation system $\Sigma$:
$$\Sigma:\cases{
{\bf x}_{i0}+{\bf x}_{i1}=a_i&(all $i<2$)\cr
{\bf x}_{0i}+{\bf x}_{1i}=b_i&(all $i<2$)\cr
}$$
It is trivial that if one of the $a_i$'s or one of the $b_j$'s
is zero, then $\Sigma$ admits a solution in $M$. Thus suppose
now that all the $a_i$'s and all the $b_i$'s are non zero.
By Lemma 1.7, there exists a [unique] conical homomorphism
$f:R^+\to M$ such that for all $i<2$, $f(\alpha_i)=a_i$ and
$f(\beta_i)=b_i$.
Let $e$ be the inclusion map from $R^+$ into
$({\bb Z}^+)^4$; it is clearly unitary. Moreover, $R$ is the
kernel of a homomorphism of abelian groups from ${\bb Z}^4$ to
${\bb Z}$ (namely $(x_0,x_1,x_2,x_3)\mapsto x_0+x_3-x_1-x_2$),
whence ${\bb Z}^4/R$ is torsion-free; thus $e$
is strongly unitary. Now let
$N=({\bb Z}^+)^4\amalg_{e,f}M$; identify $M$ with its natural image
into $N$, and let $\bar f:({\bb Z}^+)^4\to N$ be the natural
homomorphism. By Lemma 1.6, $N$ is conical and also a strong
unitary extension of $M$. Furthermore, if $(e_i)_{i<4}$ is the
natural basis of ${\bb Z}^4$, then the following array is a
refinement matrix with entries in $N$:
$$\vbox{\offinterlineskip
\halign{\tv\cc{$#$}&\tv\cc{$#$}&\tv\cc{$#$}&\tv#\cr
\omit&\omit\hrulefill&\omit\hrulefill&\omit\cr
\omit\tvi&b_0&b_1&\cr
\noalign{\hrule}
a_0&\bar f(e_0)&\bar f(e_1)&\cr
\noalign{\hrule}
a_1&\bar f(e_2)&\bar f(e_3)&\cr
\noalign{\hrule}
}}$$
Thus $\Sigma$ admits a solution in $N$.
Since $N$ is a [conical] strong unitary extension of $M$ and by
assumption on $M$, $\Sigma$ also admits a solution in $M$. Thus
$M$ is a refinement cone.\smallskip

Let us now prove (b). It will also use some unitary extension
of $M$, this time constructed directly and not by using some
amalgamated sum. The case where $a=0$ is trivial, thus let us
suppose from now on that $a\ne0$. For every real number $r$,
denote by $\uip{r}$ the least integer $\geq r$. Then equip
${\cal N}=M\times{\bb Z}^+$ with its natural (product) monoid
structure, and let $\sim$ the binary relation defined on
$\cal N$ by
$$(x,m)\sim(y,n)\Longleftrightarrow\bigl(
m\equiv n\pmod p\quad\and\quad
x+\uip{m/p}a=y+\uip{n/p}a\bigr).$$
It is easy to verify that $\sim$ is an equivalence relation
on $\cal N$, which is moreover {\it compatible with the
addition} of $\cal N$, {\it i.e.}, a {\it congruence} on
$\cal N$. Put $N={\cal N}/\!\!\sim$, and for all $(x,n)\in\cal N$,
denote by $[x,n]$ the equivalence class of $(x,n)$ modulo
$\sim$. {\it Since $a$ is non zero} and $M$ is conical, it is
easy to verify that $N$ is conical. 
Then let $j:M\to N,\ x\mapsto[x,0]$. It is easy to
verify that $j$ is a one-to-one homomorphism from $M$ into $N$.
It is unitary: indeed, it is trivially cofinal; further,
let $x,y\in M$ and let $(z,n)\in{\cal N}$
such that $j(x)+[z,n]=j(y)$; this means that
$(x+z,n)\sim(y,0)$. Thus $n=kp$ for some $k\in{\bb Z}^+$, and thus
$[z,n]=[z+ka,0]=j(z+ka)\in j[M]$.

Now put $u=[0,1]$. Then $pu=[0,p]=[a,0]=j(a)$. Furthermore,
for all $x$ and $y$ in $M$, $j(x)+u=j(y)+u$ if and
only if $(x,1)\sim(y,1)$, {\it i.e.},
$x+\uip{1/p}a=y+\uip{1/p}a$, {\it i.e.}, $x+a=y+a$.
Identify $M$ with $j[M]$. 
Then the following [finite] equation system with unknown
{\bf z}

$$\cases{
p{\bf z}=a&\cr
x+{\bf z}=y+{\bf z}&(all $x,y\in X$ such that $x+a=y+a$)\cr
}$$
admits a solution in $N$, namely $u$. Thus it admits a
solution in $M$, which completes the proof of (b).\qed\medskip

In particular, {\it every cone embeds strongly unitarily into
a refinement cone} and {\it every cone embeds unitarily into
a normally divisible refinement cone}. Now, the interest of
having considered unitary embeddings appears in the following
theorem, where we only give a sample of the kind of
results that similar methods can yield:\medskip

\Theorem1.9. Let {\bf C} be either the
class of all cancellative cones, or the class of all
separative cones. Then any {\bf C}-strongly unitarily closed
element of {\bf C} is a refinement cone; if in
addition it is {\bf C}-unitarily closed, then it is
normally divisible.\smallskip
{\rm Note that it is trivial that both classes above are closed
under unions of chains, thus Proposition 1.5 can be
used for them.}
\Proof. In each case, let $M$ be a
{\bf C}-unitarily closed (resp. {\bf C}-strongly unitarily
closed) element of {\bf C}. By Proposition 1.5 (applied to the
class of all cones), $M$ admits a unitary (resp. strong
unitary) embedding into a unitarily closed (resp. strongly
unitarily closed) cone $N$. There is no reason for $N$ to
belong to {\bf C}, thus the necessity to prove the
following\medskip

\noindent{\bf Claim.} {\sl There exists a least (for
the inclusion) congruence $\equiv$ on $N$ such
$\hat N=N/\!\!\equiv$ belongs to {\bf C}, and then the natural
homomorphism $M\to\hat N$ is unitary (resp. strongly
unitary).}\smallskip

\noindent{\bf Proof of Claim.} In the case where {\bf C} is
the class of all cancellative cones, let $\equiv$ be defined
by
$$x\equiv y\Leftrightarrow(\exists z)(x+z=y+z),$$
and in the case where {\bf C} is the class of all separative
cones, let $\equiv$ be defined by
$$x\equiv y\Leftrightarrow(\exists n\in{\bb N})
(nx+y=(n+1)x\ \and\ x+ny=(n+1)y).$$
It is well-known that in both cases, $\equiv$ is the least
congruence on $N$ with cancellative (resp. separative)
quotient monoid (details about this last fact can be found in
[{\bf\clpr}, vol. 1]).

Thus to conclude the proof of the
first part of the statement of the Claim, it suffices to prove
that in both cases, $\hat N$ is conical. For all $y\in N$,
denote by $[y]$ the equivalence class of $y$ modulo $\equiv$.
Since $N$ is conical, it suffices to prove that in both cases,
for all $y\in N$, $[y]=[0]$ implies that $y=0$. In the case
where {\bf C} is the class of all cancellative cones, there
exists $x\in N$ such that $y+x=x$; since $M$ is cofinal
in $N$, one may take $x\in M$. Since $N$ is a unitary extension
of $M$, we get $y\in M$. Since $M$ is cancellative,
$y=0$. In the case where {\bf C} is the class of all
separative cones, there exists $n\in{\bb N}$ such that
$ny+0=(n+1)y$ and $y+n0=(n+1)0$, whence $y=0$. Hence in both
cases, $\hat N$ is conical.

Now, since $M\in{\bf C}$, it is easy to see that $j$ is
one-to-one and moreover, since $M$ is cofinal in $N$, the
image of $M$ under the natural homomorphism $j:M\to\hat N$ is
also cofinal. Moreover, it is not difficult to prove both
following statements: \smallskip

\item{---} For all $x_0,x_1\in M$ and $y\in N$, if
$x_0+y\equiv x_1$, then $y\in M$.

\item{---} For all $m\in{\bb N}$, all $x\in M$ and all $y\in N$,
$my\equiv x$ implies that $y\in M$.\smallskip

(Note that in the case where {\bf C} is the class of all
cancellative cones, one uses the fact that $M$ is {\it cofinal}
in $N$).

But these two statements (together with the fact that $j$ is
one-to-one and cofinal) imply clearly unitarity (resp. strong
unitarity) of $j$.\qed\medskip

Once this is proved, the conclusion is easily reached: for
example for the finite refinement property, if $a_i$, $b_i$
($i<2$) are elements of $M$ such that $a_0+a_1=b_0+b_1$, then,
since $N$ is a refinement monoid (this results from Theorem
1.8), the equation system
$$\cases{
{\bf x}_{i0}+{\bf x}_{i1}=a_i&(all $i<2$)\cr
{\bf x}_{0i}+{\bf x}_{1i}=b_i&(all $i<2$)\cr
}$$
admits a solution in $N$, thus also in $\hat N$. Since the
natural map $M\to\hat N$ is strongly unitary, since both $M$
and $\hat N$ belong to {\bf C} and by assumption on $M$, the
equation system also admits a solution in $M$, which proves
that $M$ is a refinement monoid. The proof for normal
divisibility runs along the same lines.\qed\medskip

Using Proposition 1.5, we immediately get the following\medskip

\Corollary1.10. Let {\bf C} be either the
class of cancellative cones or the class of separative cones.
Then every element of {\bf C} admits a unitary (resp. strong
unitary) embedding into a normally divisible refinement cone
(resp. a refinement cone).\qed\medskip

\Remark1.11. Cancellative cones are of course
exactly the positive cones of \poag s. Thus, by Corollary
1.10, we get immediately the folklore result that {\it every
directed \poag\ $A$ embeds into a Riesz group $B$} (a {\it
Riesz group} is a directed interpolation group). While
unitarity gives no additional information on the embedding in
the case of groups, it allows to have $B$ {\it divisible}
(and even with divisible positive cone); on
the other hand, strong unitarity yields $B/A$ {\it
torsion-free} (but one may no longer be able to maintain
divisibility of $B$), in which case $B$ inherits whatever
torsion-freeness $A$ enjoys. In particular, {\it every
directed torsion-free \poag\ $A$ embeds into a
torsion-free Riesz group $B$ such that $B/A$ is torsion-free};
{\it see} also Corollary 1.15.

In fact, it is not difficult to prove that every 
(not necessarily directed) \poag\ $A$
embeds cofinally into an interpolation group $B$ whose
underlying group is of the form $A\oplus F$ where $F$ is free
abelian (on $|A^+|$ generators), thus yielding both results
above; {\it see} also Remark 3.13.\medskip

As an extension of Corollary 1.10 to further classes of cones,
note for example the following\medskip

\Corollary1.12. Every
stably finite separative cone embeds unitarily (resp.
strongly unitarily) into a normally divisible stably
finite separative refinement cone (resp. a stably finite
separative refinement cone).
\Proof. Note that if $N$ is a unitary extension
of $M$ and $M$ is stably finite, then $N$ is also stably
finite. Then apply Corollary 1.10.\qed\medskip

Nevertheless, as we shall now see, it is not always possible
to preserve unitarity in all these embedding theorems.
Theorem 1.14 below gives a sample of this kind of
situation.\medskip

\Definition1.13. Let
${\bb P}$ be a multiplicative subsemigroup of ${\bb N}$. Say that a
\cm\ is {\it ${\bb P}$-torsion-free} (resp.
{\it${\bb P}$-unperforated}) when for all $p\in{\bb P}$, it satisfies
the axiom $(\forall x,y)(px=py\Rightarrow x=y)$ (resp.
$(\forall x,y)(px\leq py\Rightarrow x\leq y)$). Denote by
${\bf C}_{\bb P}$ (resp. ${\bf C}'_{\bb P}$) the class of
all ${\bb P}$-torsion-free (resp. ${\bb P}$-torsion-free and
${\bb P}$-unperforated) cones.\medskip

\Theorem1.14.
\item{{\rm(a)}} Every ${\bf C}_{\bb P}$-strongly unitarily closed
cone is a refinement cone.
\item{{\rm(b)}} Every ${\bf C}'_{\bb P}$-closed cone is a normally
divisible refinement cone.
\Proof. (a) Let $M$ be a ${\bf C}_{\bb P}$-strongly unitarily closed
cone. By Proposition 1.5, $M$ embeds strongly unitarily into a 
strongly unitarily closed cone $N$. By Theorem 1.8, $N$ is a
refinement cone. Let $\equiv$ be the
congruence on $N$ defined by
$$x\equiv y\Leftrightarrow(\exists p\in{\bb P})(px=py).$$
It is trivial that $N/\!\!\equiv$ is conical.
By definition, $N/\!\!\equiv$ is ${\bb P}$-torsion-free.
It is also easy to see that the natural homomorphism
$M\to N/\!\!\equiv$ is strongly unitary.

Now the argument used (for the refinement) at the end of the
proof of Theorem 1.9 allows us easily to conclude that $M$ is
a refinement cone.\smallskip

Let us now prove (b). Thus let $M$ be a ${\bf C}'_{\bb P}$-closed
cone. By Proposition 1.5, $M$ order-embeds into a closed cone
$N$. By Theorem 1.8, $N$ is a normally divisible
refinement cone. As above, let $\equiv$ be the
congruence on $N$ defined by
$$x\equiv y\Leftrightarrow(\exists p\in{\bb P})(px=py).$$
It is again trivial that $N/\!\!\equiv$ is conical.
By definition, $N/\!\!\equiv$ is ${\bb P}$-torsion-free; {\it by
using the fact that $N$ is divisible}, one can also see easily
that $N/\!\!\equiv$ is ${\bb P}$-unperforated. Finally, one checks
easily that the natural homomorphism $M\to N/\!\!\equiv$ is an
order-embedding. From then on, the usual argument shows that
$M$ is a normally divisible refinement monoid.
\qed\medskip

\Corollary1.15.
\item{{\rm(a)}} Every ${\bb P}$-torsion-free
cone embeds strongly unitarily into a ${\bb P}$-torsion-free
refinement cone.
\item{{\rm(b)}} Every directed ${\bb P}$-torsion-free \poag\ $A$
embeds cofinally into a ${\bb P}$-torsion-free Riesz group $B$
such that $B/A$ is torsion-free.
\item{{\rm(c)}} Every ${\bb P}$-torsion-free
${\bb P}$-un\-per\-fo\-ra\-ted cone order-embeds into a
${\bb P}$-torsion-free ${\bb P}$-un\-per\-fo\-ra\-ted normally
divisible refinement cone.
\item{{\rm(d)}} Every directed ${\bb P}$-un\-per\-fo\-ra\-ted
\poag\ $A$ embeds cofinally into a ${\bb P}$-un\-per\-fo\-ra\-ted
Riesz group $B$ such that $B/A$ is torsion-free.
\Proof. (a), (c) result immediately from
Proposition 1.5 and Theorem 1.14, while (b) results
immediately from Remark 1.11. Let us see now (d); thus let
{\bf C} be the class of all ${\bb P}$-unperforated cancellative
cones and let $M$ be a {\bf C}-strongly unitarily closed cone;
thus $M=A^+$ for some directed \poag\ $A$. By (b) above, $A$
embeds cofinally into a ${\bb P}$-torsion free Riesz group $B$
such that $B/A$ is torsion-free. Now define a subset $P$ of
$B$ by
$$P=\{x\in B:\ (\exists p\in{\bb P})(px\geq0)\}.$$
Using the fact that $B$ is ${\bb P}$-torsion-free, one sees
easily that $P$ is the positive cone of a ${\bb P}$-unperforated
\poag\ $B'$ of underlying group $B$; in particular,
$P\in{\bf C}$. Since $A$ is ${\bb P}$-unperforated, the inclusion
map from $A$ into $B'$ is an order-embedding;
since $P$ contains $B^+$, it is
still cofinal; thus, since $B/A$ is torsion-free,
the inclusion map from $M$ into $P$ is strongly unitary. If
$a_i$, $b_i$ ($i<2$) are elements of $M$ such that
$a_0+a_1=b_0+b_1$, then, since $B^+$ is a refinement cone, the
equation system
$$\Sigma:\cases{
{\bf x}_{i0}+{\bf x}_{i1}=a_i&(all $i<2$)\cr
{\bf x}_{0i}+{\bf x}_{1i}=b_i&(all $i<2$)\cr
}$$
admits a solution in $B^+$, thus {\it a fortiori} in $P$. By
assumption on $M$, $\Sigma$ admits a solution in $M$; whence
$M$ is a refinement cone. We conclude again by
Proposition 1.5.\qed\medskip

Note that one cannot improve ``refinement cone" into
``divisible refinement cone" in (a) above (because every
divisible torsion-free \cm\ is unperforated while there are
torsion-free perforated cones); however, this drawback will be
overcome in section 3 where we will introduce {\it
quasi-divisibility}. Note also that one cannot improve
``order-embeds" into ``embeds unitarily" in (c) above: for
example, $M={\bb Z}^+\cup\{\infty\}$ does not embed unitarily
into any divisible, torsion-free unperforated \cm\ (because in
such an extension, multiplication by positive rational numbers
would be defined, and then $(1/2)+\infty=\infty$ while
$1/2\notin M$).\smallskip 

We shall finally discuss briefly another property, the axiom
(WSD) (``Weak Sum Decomposition"), slightly weaker than the
axiom (SD) considered in [{\bf\wehrIII}]:
$$\displaylines{
({\rm WSD}):\qquad(\forall a_0,a_1,b,c)\bigl[a_0+a_1+c=b+c
\Rightarrow\hfill\cr
\hfill(\exists x_0,x_1)(a_0+c=x_0+c\ \and\ a_1+c=x_1+c\ \and\
b=x_0+x_1)\bigr].\qquad\cr
}$$

It is easy to verify that {\it every cone
satisfying {\rm (WSD)} is antisymmetric} ({\it i.e.}, it
satisfies
$(\forall x,y)((x\leq y\ \and\ y\leq x)\Rightarrow x=y)$),
but it should be noted that (WSD) is {\it not} a
consequence of the finite refinement property
plus antisymmetry. A simple
counterexample for this is the monoid 
${\bf\Lambda}({\bb Q}^+)$ of all nonempty lower subsets of
${\bb Q}^+$:
for any non-negative real number $r$, identify $r$ with the
interval $[0,\ r]\cap{\bb Q}^+$, and put $r^-=[0,\ r)\cap{\bb Q}^+$ (so
that for irrational $r$, we have $r=r^-$); then, let $\alpha$
be any irrational number such that $0<\alpha<1$, and take
${\goth a}_0=\alpha$, ${\goth a}_1=1-\alpha$, ${\goth b}=1$ and
${\goth c}=1^-$; then ${\goth a}_0$, ${\goth a}_1$, ${\goth
b}$, ${\goth c}$ witness failure of (WSD) in ${\bf\Lambda}({\bb Q}^+)$;
{\it see} also the more general [{\bf\wehrIII}, Theorem 2.21].
Another source of counterexamples (that refinement does not
imply (WSD)) comes from the fact that
refinement plus (WSD) implies that the maximal cancellative
quotient has refinement (it implies in fact (HREF) as defined
in [{\bf\wehrIII}], but it is not equivalent to (HREF); for
example, ${\bf\Lambda}({\bb Q}^+)$ satisfies (HREF) [{\bf\wehrIII},
Theorem 2.11]): so for J. Moncasi's example of a regular
ring $R$ such that $K_0(R)$ is not a Riesz group, the cone
of isomorphism classes of finitely generated projective
right $R$-modules satisfies refinement but not (WSD);
similarly, if $A$ and $B$ are Riesz groups such that
$A\otimes B$ is not an interpolation group [{\bf\wehrIV},
Example 1.4], the tensor product of $A^+$ and $B^+$ as \cm s
satisfies refinement [{\bf\wehrIV}, Theorem 2.9], but not
(WSD). Our next proposition will produce cones satisfying
(WSD).

\medskip

\Proposition1.16. Let {\bf C} be the
class of all antisymmetric cones. Then every
{\bf C}-closed cone satisfies {\rm(WSD)}.
\Proof. Let $M$ be a {\bf C}-closed
cone. Let $a_0$, $a_1$, $b$ and $c$ in $M$
such that $a_0+a_1+c=b+c$. Consider the following
equation system with unknowns ${\bf x}_0$ and ${\bf x}_1$:
$$\Sigma:\cases{
{\bf x}_0+c=a_0+c;\cr
{\bf x}_1+c=a_1+c;\cr
{\bf x}_0+{\bf x}_1=b.\cr
}$$
We shall construct an extension
$N$ of $M$ in {\bf C} where $\Sigma$ admits a solution.

If $a_0=0$ then $(x_0=0;x_1=b)$ is a solution of $\Sigma$ in
$M$. Similarly for $a_1=0$. If $b=0$, then, since $M$ is
antisymmetric, $a_0+c=a_1+c=c$, thus $x_0=x_1=0$ is a solution
of $\Sigma$ in $M$.

So now, suppose that $a_0$, $a_1$ and $b$ are non zero. Let
$F={\bb Z}^+\times{\bb Z}^+$ and let $f:F\to M$ be the unique
homomorphism of monoids sending $e_0=(1,0)$ on $a_0$ and
$e_1=(0,1)$ on $a_1$. Let $\to$ the binary relation
(``rewriting rule") on ${\cal N}=M\times F$ consisting exactly
on all pairs of the form
$$\displaylines{ (x,r)\to(x,r),\cr
(x,e_i+r)\to(x+a_i,r)\qquad({\rm for}\ i<2)\qquad{\rm if}\
x\geq c,\cr (x,e_0+e_1+r)\to(x+b,r).\cr }$$

It is immediate to verify that $\to$ is compatible with the
addition of $\cal N$ ({\it i.e.}, $\xi\to\eta$ implies
$\xi+\zeta\to\eta+\zeta$).\medskip

\noindent{\bf Claim.} {\sl The relation $\to$ is {\rm
confluent}, {\rm i.e.}, for all $\xi$, $\eta_0$, $\eta_1$ in
$\cal N$ such that $\xi\to\eta_i$ (all $i<2$), there exists
$\zeta\in\cal N$ such that $\eta_i\to\zeta$ (all
$i<2$).}\smallskip

\noindent{\bf Proof of Claim.} There are essentially two cases
that are not completely trivial to consider:\smallskip

\item{\bfi Case 1.} {\it $\xi=(x,r)$ with $r=e_0+r_0=e_1+r_1$
and $x\geq c$, and $\eta_i=(x+a_i,r_i)$ (all $i<2$).}
Then there exists $r'\in F$ such that $r_i=e_{1-i}+r'$ (all
$i<2$). Take $\zeta=(x+a_0+a_1,r')$.
\smallskip

\item{\bfi Case 2.} {\it$\xi=(x,e_0+e_1+r)$ with $x\geq c$,
$\eta_0=(x+a_0,e_1+r)$ and $\eta_1=(x+b,r)$.}
Then since $x\geq c$, we have $x+a_0+a_1=x+b$, and then it is
easy to verify that one can take
$\zeta=\eta_1=(x+b,r)$.\qed\smallskip

Now let $\to^*$ be the transitive closure of $\to$. Since
$\to$ is reflexive, $\to^*$ is reflexive; since
$\to$ is compatible with the addition, $\to^*$ is
also compatible with the addition, and since
$\to$ is confluent, $\to^*$ is also confluent. Thus the
binary relation $\equiv$ defined on $\cal N$ by
$$\xi\equiv\eta\Leftrightarrow(\exists\zeta)
(\xi\to^*\zeta\ \and\ \eta\to^*\zeta)$$
is a congruence on $\cal N$. Let $N={\cal N}/\!\!\equiv$, and
for all $(x,r)\in\cal N$, denote by $[x,r]$ its equivalence
class modulo $\equiv$. To prove that $N$ is conical, it
suffices to prove that for all $(x,r)\in M\times F$,
$[x,r]=[0,0]$ implies that $(x,r)=(0,0)$. Thus suppose
$[x,r]=[0,0]$. By definition, there exists
$(y,s)\in M\times F$ such that $(x,r)\to^*(y,s)$ and
$(0,0)\to^*(y,s)$. The second condition easily implies that
$(y,s)=(0,0)$, whence $(x,r)\to^*(0,0)$. Since $\to$ increases
the first coordinate, $x\leq0$, thus, by conicality of $M$,
$x=0$, so that $(0,r)\to^*(0,0)$. But since $a_0$, $a_1$ and
$b$ are non zero and $M$ is conical, this is possible only when
$r=0$; whence $N$ is conical.

Further, let $j:M\to N,\ x\mapsto[x,0]$.
It is immediate that $j$ is a monoid homomorphism. If $x$ and
$y$ are two elements of $M$ such that $j(x)=j(y)$, then there
exists $(z,r)\in\cal N$ such that $(x,0)\to^*(z,r)$ and
$(y,0)\to^*(z,r)$. Since $\to$ strictly decreases the second
coordinate (except in the cases where it is equality), we
necessarily have $r=0$ and $z=x=y$. Now if we just suppose
that $j(x)\leq j(y)$, then there exists $(z,r)\in\cal N$ such
that $(x+z,r)\equiv(y,0)$, thus there exists
$(y',r')\in\cal N$ such that $(x+z,r)\to^*(y',r')$ and
$(y,0)\to^*(y',r')$. Thus, as before, $r'=0$ and $y'=y$; since
$x+z\leq y'$, we deduce $x\leq y$. Therefore, $j$ is an
order-embedding from $M$ into $N$.

Further, put $x_i=[0,e_i]$ for all $i<2$. Then
$x_i+j(c)=[c,e_i]=[c+a_i,0]=j(a_i)+j(c)$. Furthermore,
$x_0+x_1=[0,e_0+e_1]=[b,0]=j(b)$.
Thus, the image of $\Sigma$ under $j$ admits a solution in
$N$.

Finally, let $N'$ be the maximal antisymmetric quotient of
$N$: that is, $N'=N/\!\!\equiv$ where $\equiv$ is defined by
$x\equiv y\Leftrightarrow(x\leq y\ \and\ y\leq x)$ (it is
obviously a congruence on $N$), let $\pi:N\to N'$ be the
natural projection. Then $\pi\circ j$ is an order-embedding
from $M$ into $N'$, and the image of $\Sigma$ under
$\pi\circ j$ admits a solution in $N'$. Since $M$ is
{\bf C}-closed, $\Sigma$ also admits a solution in $M$, which
concludes the proof.\qed\medskip

\Corollary1.17. Every antisymmetric cone
order-embeds into a normally divisible refinement cone
satisfying {\rm(WSD)}.\qed\medskip

Note that it is easy to verify that every normally
divisible refinement monoid satisfying (WSD) satisfies in fact
the axiom (SD) of [{\bf\wehrIII}]: thus it is a {\it
refinement algebra} in the sense of [{\bf\wehrIII}].
\bigskip\goodbreak

\sectionheadline{\S2. Order-embedding simple cones.}

In this section we will try to extend the results of previous
section to {\it simple} cones. A drawback to this is that we
will need to use {\it divisibility} in several cases, even for
the mere proof of the finite refinement property; thus this
section will contain no result concerning strongly unitary
embeddings. We will return back to the latter in section 3, to
solve the remaining questions. We shall make use of the
following lemma [{\bf\wehr}, Lemma 1.9]:\medskip

\Lemma2.1. Let $n\in{\bb N}$, let $M$ be a
refinement monoid and let $a$, $b$ and $c$ be elements of $M$
such that $a+b=nc$. Then there exist elements $c_k$
($0\leq k\leq n$) of $M$ such that
$$a=\sum_{k\leq n}kc_k,\quad b=\sum_{k\leq n}(n-k)c_k\quad
\and\quad c=\sum_{k\leq n}c_k.\qeddis\medskip

For all elements $a$ and $b$ of an algebraic \cm\ $M$, we will
write $a\propto b$ when there exists $n\in{\bb N}$ such that
$a\leq nb$, and $a\asymp b$ when $a\propto b$ and $b\propto a$.
Furthermore, for all $a\in M$, put
$M(a)=\{x\in M:\ x\asymp a\}\cup\{0\}$.
A subset $X$ of $M$ is {\it$\asymp$-trivial} when
$(\forall x,y\in X)(x\asymp y)$. Say that a refinement
matrix is $\asymp$-trivial when the set of all its entries
is $\asymp$-trivial. 
\medskip

\Lemma2.2. Let $M$ be a refinement monoid.
Then every single $\asymp$-equivalence class of $M$ is
downward directed.
\Proof. It suffices to prove that if $a$ and $b$
are two elements of $M$ such that $a\asymp b$, then there
exists $c\asymp a$ such that $c\leq a$ and $c\leq b$. By
assumption there exists $n\in{\bb N}$ such that $a\leq nb$. By
Lemma 2.1, there are elements $c_k$ ($0\leq k\leq n$) such that
$a=\sum_{k\leq n}kc_k$ and $b=\sum_{k\leq n}c_k$. Put
$c=\sum_{k=1}^nc_k$. Then $c\leq a$ and $c\leq b$; thus
$c\propto a$, but $a\leq nc$, thus $a\asymp c$.
\qed\medskip

We now come to the main lemma of this section:\medskip

\Lemma 2.3. Let $M$ be a cone and let
$a\in M$. Then $M(a)$ is a simple cone. If in addition $M$ is
a normally divisible refinement cone, then $M(a)$ is a simple
normally divisible refinement cone.\smallskip
{\rm Note that in fact, in order to get refinement in $M(a)$,
the hypothesis of the second paragraph of Lemma 2.3 can be
weakened into ``$M$ is a refinement cone satisfying the axiom
$$(\forall a,b,c)\bigl(a+c=b+c\Rightarrow(\exists x)
(2x=c\ \and\ a+x=b+x)\bigr)",$$
as the proof will show it.}
\Proof. It is obvious that $M(a)$ is a submonoid
of $M$. Now let $x$ and $y$ be two non zero elements of
$M(a)$. Since $x\asymp y$, there are $n\in{\bb N}$ and $z\in M$
such that $x+z=ny$. Therefore, $y+z\in M(a)$ and
$x+(y+z)=(n+1)y$, so that $M(a)$ satisfies $x\propto y$;
therefore, $y$ is an order-unit of $M(a)$; whence $M(a)$ is
simple.\smallskip

From now on suppose that $M$ is a normally divisible
refinement monoid. Let $a_0$, $a_1$, $b_0$ and $b_1$ be
elements of $M(a)$ such that $a_0+a_1=b_0+b_1$; we must find a
refinement of this equality in $M(a)$. If one of the
$a_i$'s or one of the $b_i$'s is equal to $0$, then the
problem obviously admits a solution, thus suppose that none
of the $a_i$, $b_i$'s is equal to zero.
Put $X=\{a_0,a_1,b_0,b_1\}$. By Lemma
2.2, there exists $c\asymp a$ such that $c\leq X$. Thus
there are elements $a'_i$, $b'_i$ ($i<2$) of $M$ such that
$a_i=c+a'_i$ and $b_i=c+b'_i$ (all $i<2$), whence
$a'_0+a'_1+2c=b'_0+b'_1+2c$. Now, two successive applications
of normal divisibility of $M$ yield two elements $u_0$ and
$u_1$ of $M$ such that $2u_0=2u_1=c$ and
$a'_0+a'_1+u_0+u_1=b'_0+b'_1+u_0+u_1$. Thus $u_0\asymp u_1$,
thus, since $M$ is a refinement monoid and by Lemma 2.2, there
exists $u\in M$ such that $u\leq\{u_0,u_1\}$ and
$u\asymp\{u_0,u_1\}$. Thus there are elements $u'_0$ and $u'_1$
of $M$ such that $u_i=u+u'_i$ (all $i<2$). Then one more
application of normal divisibility yields an element $v$ of $M$
such that $2v=u$. Now we have the following refinement matrix
$$\vcenter{\offinterlineskip
\halign{\tv\cc{$#$}&\tv\cc{$#$}&\tv\cc{$#$}&\tv#\cr
\omit&\omit\hrulefill&\omit\hrulefill&\omit\cr
\omit\tvi&u_0&u_1&\cr \noalign{\hrule}
u_0&v+u'_0&v&\cr
\noalign{\hrule}
u_1&v&v+u'_1&\cr
\noalign{\hrule}
}}$$
and it is $\asymp$-trivial (because
$v\asymp u\asymp u_0\asymp u_1$).
Furthermore, since $M$ is a refinement monoid, there exists a
refinement matrix with entries in $M$ of the following
form:
$$\vcenter{\offinterlineskip
\halign{\tv\cc{$#$}&\tv\cc{$#$}&\tv\cc{$#$}&\tv#\cr
\omit&\omit\hrulefill&\omit\hrulefill&\omit\cr
\omit\tvi&b'_0+u_0&b'_1+u_1&\cr
\noalign{\hrule}
a'_0+u_0&c'_{00}&c'_{01}&\cr
\noalign{\hrule}
a'_1+u_1&c'_{10}&c'_{11}&\cr
\noalign{\hrule}
}}$$
Therefore, we have the following refinement matrix in $M$:
$$\vcenter{\offinterlineskip
\halign{\tv\cc{$#$}&\tv\cc{$#$}&\tv\cc{$#$}&\tv#\cr
\omit&\omit\hrulefill&\omit\hrulefill&\omit\cr
\omit\tvi&b_0&b_1&\cr
\noalign{\hrule}
a_0&c'_{00}+v+u'_0&c'_{01}+v&\cr
\noalign{\hrule}
a_1&c'_{10}+v&c'_{11}+v+u'_1&\cr
\noalign{\hrule}
}}$$
and it is $\asymp$-trivial (all its entries
are $\geq v$); thus all its entries belong to $M(a)$.\smallskip

Let us finally verify normal divisibility. Thus let $p\in{\bb N}$
and let $b\in M(a)$ and $X\subseteq M(a)$ be a finite subset.
Put $Y=\{(x,y)\in X\times X:\ x+b=y+b\}$. Since $M$ is normally
divisible, there exists $c\in M$ such that $pc=b$ and
$(\forall (x,y)\in Y)(x+c=y+c)$. But $c\asymp b$, thus
$c\in M(a)$, and we are done.
\qed\medskip

\Remark2.4. It is not difficult to prove that
{\it for every simple atomless refinement cone $M$, the
semigroup $M^{>0}=\{x\in M:\ x>0\}$ satisfies the finite
refinement property} (this fails of course for $M={\bb Z}^+$, which
is atomic); in fact, the proof works for atomless refinement
cones $M$ that are in addition {\it prime}, {\it i.e.},
$M^{>0}$ is downward directed. A proof of this can be found in
Bergman's notes [{\bf\berg}]. This can also be considered as
a vindication of the result of Lemma 2.3.\medskip

Before going on, let us rephrase for convenience the version
of Proposition 1.5 for simple \cm s:\medskip

\Proposition2.5. Let {\bf C} be a class of
\cm s closed under unions of chains. Then every simple \cm\ 
of {\bf C} admits an order-embedding (resp. a unitary
embedding, a strong unitary embedding) into a simply
{\bf C}-closed (resp. a simply {\bf C}-unitarily closed, a
simply {\bf C}-strongly unitarily closed) element of
{\bf C}.
\Proof. Just apply Proposition 1.5 to the class
of all simple elements of {\bf C} (the point is that the union
of any chain of simple \cm s is a simple \cm).\qed\medskip

Note that the class of all cones (resp. of all cancellative,
of all separative cones) is closed under unions of chains, thus
satisfies the hypothesis of Proposition 2.5.\smallskip

Now we can prove our main embedding theorem
for simple refinement cones:\medskip

\Theorem2.6. Let {\bf C} be the class
of all cones (resp. {\rm cancellative}, resp. {\rm separative}
cones). Then every simply {\bf C}-unitarily
closed element of {\bf C} is a normally divisible refinement
cone.
\Proof. Let $M$ be a simply {\bf C}-unitarily
closed element of {\bf C}. By Proposition 1.5,
$M$ admits a unitary embedding into a {\bf C}-unitarily
closed element $N$ of {\bf C}. By
Theorems 1.8 and 1.9, $N$ is a normally divisible refinement
cone.\smallskip

Let us first prove
that $M$ is a refinement monoid. Thus let $a_i$, $b_i$ ($i<2$)
be elements of $M$ such that $a_0+a_1=b_0+b_1$. As usual, let
$\Sigma$ be the following equation system with unknowns
${\bf x}_{ij}$ ($i,j<2$):
$$\Sigma:\cases{
{\bf x}_{i0}+{\bf x}_{i1}=a_i&(all $i<2$)\cr
{\bf x}_{0i}+{\bf x}_{1i}=b_i&(all $i<2$)\cr
}$$
If one of the $a_i$, $b_i$'s is equal to zero,
then it is trivial that $\Sigma$ admits a solution in $M$.

So now, suppose that neither of the $a_i$, $b_i$'s is equal to
zero. Then, since $M$ is simple, all four of them are
order-units, thus $a_0\asymp a_1\asymp b_0\asymp b_1$. By
Lemma 2.3, $\Sigma$ admits a solution in $N(a_0)$. Since $M$
is simple, we have $M\subseteq N(a_0)$; since $N$ is a unitary
extension of $M$, so is $N(a_0)$. Since $N(a_0)\in{\bf C}$ and
$M$ is simply {\bf C}-unitarily closed, $\Sigma$ also admits a
solution in $M$. Thus $M$ is a refinement monoid.\smallskip

Let now $X$ be a finite subset of $M$ and let $a\in M$.
Put $Y=\{(x,y)\in X\times X:\ x+a=y+a\}$ and consider the
following equation system with unknown {\bf z}
$$\Sigma':\cases{
p{\bf z}=a&\cr
x+{\bf z}=y+{\bf z}&(all $(x,y)\in Y$)\cr
}$$
If $a=0$, then it is obvious that $\Sigma'$ admits a solution
in $M$, thus suppose that $a\ne0$. Thus all
parameters of $\Sigma'$ belong to $N(a)$, thus, since $N(a)$
is normally divisible (by Lemma 2.3), $\Sigma'$ admits a
solution in $N(a)$. But as for refinement, $N(a)$ is a unitary
extension of $M$ which belongs to {\bf C}, whence $\Sigma'$
admits a solution in $M$. Hence $M$ is normally
divisible.\qed\medskip

\Corollary2.7. Let {\bf C} be the class
of all cones (resp. {\rm cancellative}, resp. {\rm separative}
cones). Then every simple element of {\bf C} embeds
unitarily into a simple normally divisible refinement
cone belonging to {\bf C}.\qed\medskip

By applying again Remark 1.11, one can easily extend this
result to other classes of cones. For example, by taking for
{\bf C} the class of all cancellative cones, one
obtains easily the following result:\medskip

\Corollary2.8. Every simple \poag\
embeds into a simple Riesz group with divisible positive
cone.\qed\medskip

(Note also that an [order-] embedding of {\it simple} \poag s
is necessarily cofinal).

In fact, it is not difficult to extend this result to the class
of unperforated \poag s:\medskip

\Corollary2.9. Every
simple unperforated \poag\ embeds into a simple
divisible dimension group.
\Proof. Denote by {\bf C} the class of all
cancellative cones ({\it i.e.}, the class of all
positive cones of \poag s), and by ${\bf C}'$ the class of all
unperforated elements of {\bf C}. By Proposition 2.5, it
suffices to prove that every simply ${\bf C}'$-closed cone $M$
is a divisible refinement cone. By Corollary 2.8, $M$
order-embeds into the positive cone $N$ of some simple
Riesz group $H$, with $N$ divisible. Now define a convex
subgroup $I$ of $H$ by
$$I=\{x\in H:\ (\exists m\in{\bb N})(mx=0)\},$$
and let $\hat N$ be the positive cone of $H/I$.
Then it is obvious, using torsion-freeness of $M$, that the
natural homomorphism $j:M\to\hat N$ is one-to-one. It is
also easy to see, {\it using unperforation}, that $j$ is an
order-embedding. Furthermore, using the fact that $N$ is
divisible, it is not difficult to infer that $\hat N$ is
unperforated (thus belongs to ${\bf C}'$) and divisible.

Now let $a\in M$ and let $m\in{\bb N}$. Then, since $N$ is
divisible, the equation $m{\bf x}=a$ admits a solution in $N$,
thus in $\hat N$; thus, by assumption on $M$, the equation
also admits a solution in $M$. Therefore, $M$ is divisible.
The proof for refinement is similar ({\it see} also the proof
of Theorem 1.14 (b)).\qed\medskip

A similar application of Corollary 2.7 yields also, with a
proof similar to the one of Theorem 1.14 (b), the following
result (recall that ${\bb P}$ is a given multiplicative
subsemigroup of ${\bb N}$):\medskip

\Corollary2.10. Every
simple ${\bb P}$-torsion-free ${\bb P}$-unperforated cone order-embeds
into a normally divisible simple ${\bb P}$-torsion-free
${\bb P}$-unperforated refinement cone.\qed\medskip

\Remark2.11. For arbitrary {\it torsion-free}
groups, the argument above does not work {\it a priori} because
the natural map from $M$ to $N/\!\!\equiv$ does not seem to be
an order-embedding. A deeper reason for the embedding problem
to be more difficult in the torsion-free case is that {\it
every divisible torsion-free \cm\ is unperforated}, while there
are perforated torsion-free cones. But normal divisibility
plays an important role in the proof of Corollary 2.7, even
for the part about refinement!
However, we will see in next chapter that the
result of Corollary 2.9 still holds for
torsion-free \poag s (and even more general classes of
\cm s) as far as refinement is concerned.\medskip

Finally, let us prove the ``simple analogue" of Proposition
1.16:\medskip

\Proposition2.12. Let {\bf C} be the
class of all antisymmetric cones. Then every simply
{\bf C}-closed cone satisfies {\rm(WSD)}.\smallskip
Therefore, {\it every simple antisymmetric cone order-embeds
into a simple normally divisible refinement cone satisfying
{\rm(WSD)}}.
\Proof. Let $M$ be a simply {\bf C}-closed
cone. By Proposition 1.5, $M$ order-embeds into a
{\bf C}-closed cone $N$. By Proposition 1.16, $N$ satisfies
(WSD), while by Theorem 2.6, $M$ is a normally divisible
refinement cone. Now let $a_0$, $a_1$, $b$ and $c$ be
elements of $M$ such that $a_0+a_1+c=b+c$. Consider the
following equation system with unknowns ${\bf x}_0$ and
${\bf x}_1$:
$$\Sigma:\cases{
{\bf x}_0+c=a_0+c;\cr
{\bf x}_1+c=a_1+c;\cr
{\bf x}_0+{\bf x}_1=b.\cr
}$$
Our goal is to prove that $\Sigma$ admits a solution in $M$.
One deals with the easy cases $a_0=0$ or $a_1=0$ or $b=0$
just as at the beginning of the proof of Proposition 1.16.

Thus suppose from now on that $a_0$, $a_1$ and $b$ are all non
zero. By Lemma 2.2 and Theorem 2.6, there exists $d\ne0$ in $M$
such that $d\leq a_i$ (all $i<2$) and $2d\leq b$. Thus there
are elements $a'_i$ ($i<2$) and $b'$ of $M$ such that
$a_i=a'_i+d$ (all $i<2$) and $b=b'+2d$. Therefore,
$a'_0+a'_1+c+2d=b'+c+2d$, whence, applying twice the normal
divisibility of $M$, there are elements $u_i$ ($i<2$) of $M$
such that $2u_i=d$ (all $i<2$) and
$a'_0+a'_1+u_0+u_1+c=b'+u_0+u_1+c$. Since $N$ satisfies (WSD),
there are elements $y_0$ and $y_1$ of $N$ such that
$a'_i+u_i+c=y_i+c$ (all $i<2$) and $b'+u_0+u_1=y_0+y_1$. Put
$x_i=y_i+u_i$. Then it is easy to verify that $(x_0,x_1)$ is a
solution of $\Sigma$ in $N$. Furthermore,
$u_i\leq x_i\leq a_i+c$ and $2u_i=d$ with $d>0$ in $M$, thus,
since $M$ is simple, $x_i\asymp d$. Thus $\Sigma$ admits a
solution in $N(d)$ which is a simple (conical) extension of
$M$, whence, by assumption on $M$, $\Sigma$ admits a solution
in $M$.\qed\medskip

\bigskip\goodbreak

\sectionheadline{\S3. The torsion-free case.}

As we have seen in Remark 2.11, the torsion-free case cannot
be handled directly by the methods of previous sections, the
main problem being divisibility. Thus we will need to make up
for this by introducing a new condition, weaker than
divisibility, which we will naturally call as in
[{\bf\wehrIII}, definition 2.19] {\it quasi-divisibility}.
\smallskip

{\it Throughout this section, we will fix a
multiplicative subsemigroup ${\bb P}$ of ${\bb N}$ such that
${\bb P}\not\subseteq\{1\}$.}\medskip

\Lemma3.1. Every ${\bb P}$-torsion-free \cm\
is separative.
\Proof. Fix an element $p$ of
${\bb P}\setminus\{1\}$. If $M$ is ${\bb P}$-torsion-free, then
for all $a,b\in M$ such that $2a=a+b=2b$, then an easy
induction proof shows that $(k+l)a=ka+lb=(k+l)b$ for all
$k,l\in{\bb N}$. In particular, $pa=pb$; whence $a=b$. Thus $M$ is
separative.\qed\medskip

\Definition3.2. A \cm\ is {\it
quasi-divisible} when it satisfies the following axiom:
$$(\forall x)(\exists u,v)(2u+3v=x).$$\medskip

\Lemma3.3. Every quasi-divisible \cm\
satisfies the statement
$$(\forall x)(\exists y)(2y\leq x\leq 3y).$$
\Proof. Let $u$, $v$ such that $x=2u+3v$.
Take $y=u+v$.\qed\medskip

\Lemma3.4. Let $M$ be a separative
quasi-divisible refinement cone, let $a\in M$. Then $M(a)$ is
a simple, separative quasi-divisible refinement cone.
\Proof. Simplicity has already been proved in
2.3. Separativity is trivial.
Now let $b\in M(a)$, we will find solutions in $M(a)$ of the
equation $2u+3v=b$. Since $M$ is quasi-divisible, three
successive applications of Lemma 3.3 yield easily an element
$c$ of $M$ such that $5c\leq b$ and $c\asymp b$, so that
$c\in M(a)$; let $b'\in M$ such that $b=5c+b'$. Since $M$ is
quasi-divisible, there exist $u'$ and $v'$ in $M$ such that
$b'=2u'+3v'$. Now note that $b=2u+3v$ where both $u=c+u'$ and
$v=c+v'$ belong to $M(a)$; thus $M(a)$ is quasi-divisible.

Let us check now that $M(a)$ is a refinement monoid.
Thus let $a_i$, $b_i$ ($i<2$) be elements of $M(a)$ such that
$a_0+a_1=b_0+b_1$; consider as usual the following equation
system
$$\Sigma:\cases{
{\bf x}_{i0}+{\bf x}_{i1}=a_i&(all $i<2$)\cr
{\bf x}_{0i}+{\bf x}_{1i}=b_i&(all $i<2$)\cr
}$$

If $a_i=0$ or $b_i=0$ for some $i$, then it is obvious that
$\Sigma$ admits a solution in $M(a)$. Thus suppose that all
the $a_i$, $b_i$'s are non zero. Thus
$a_0\asymp a_1\asymp b_0\asymp b_1$.
Let $X=\{a_0,a_1,b_0,b_1\}$.
By Lemma 2.2, there exists $c\asymp a$ such that $c\leq X$.
Since $M$ is quasi-divisible and by two successive applications
of Lemma 3.3, there exists $d\in M$ such that $3d\leq c$ and
$d\asymp c$. Thus there exist elements $a'_i$, $b'_i$ ($i<2$)
of $M$ such that $a_i=3d+a'_i$ and $b_i=3d+b'_i$. Hence
$a'_0+a'_1+6d=b'_0+b'_1+6d$. {\it Now, we use the separativity
of $M$} ({\it see} Lemma 3.1): this allows us to infer that
$a'_0+a'_1+d=b'_0+b'_1+d$. Therefore, one can form a
refinement matrix with entries in $M$ as follows:
$$\vbox{\offinterlineskip
\halign{\tv\cc{$#$}&\tv\cc{$#$}&\tv\cc{$#$}&\tv#\cr
\omit&\omit\hrulefill&\omit\hrulefill&\omit\cr
\omit\tvi&b'_0+d&b'_1&\cr
\noalign{\hrule}
a'_0+d&p&q&\cr
\noalign{\hrule}
a'_1&r&s&\cr
\noalign{\hrule}
}}$$
Thus the following is also a refinement matrix:
$$\vbox{\offinterlineskip
\halign{\tv\cc{$#$}&\tv\cc{$#$}&\tv\cc{$#$}&\tv#\cr
\omit&\omit\hrulefill&\omit\hrulefill&\omit\cr
\omit\tvi&b_0&b_1&\cr
\noalign{\hrule}
a_0&p+d&q+d&\cr
\noalign{\hrule}
a_1&r+d&s+2d&\cr
\noalign{\hrule}
}}$$
and all its entries belong to $M(a)$; whence $M(a)$ is a
refinement monoid.\qed\medskip

Now we come to the main theorem of this section.\medskip

\Theorem3.5. Let $M$ be
${\bf C}_{\bb P}$-strongly unitarily closed. Then $M$ is a
quasi-divisible refinement cone.\smallskip
{\rm Recall (Lemma 3.1) that $M$ is also} separative.
\Proof. We shall use throughout this proof the
techniques of amalgamated sums used in section 1.\smallskip

Let us prove quasi-divisibility.
Thus let $a\in M$; consider the equation
$\Sigma:\ 2{\bf x}+3{\bf y}=a$. If $a=0$, then $\Sigma$
trivially admits a solution in $M$. Thus suppose that $a\ne0$.
Let $D$ be the subgroup of ${\bb Z}^2$
generated by $(2,3)$, let $e:D^+\to({\bb Z}^+)^2$ be the inclusion
map. Let $f:D^+\to M$ be the unique homomorphism such that
$f\bigl((2,3)\bigr)=a$; note that $f$ is conical.
Let $N=({\bb Z}^+)^2\amalg_{e,f}M$, let
$\bar e:M\to N$ and $\bar f:({\bb Z}^+)^2\to N$ be the natural
homomorphisms. Since the map $n\mapsto(n,n)+D$ is easily seen
to be an isomorphism from ${\bb Z}$ onto ${\bb Z}^2/D$, the latter is
torsion-free; it follows easily that $e$ is strongly unitary.
By Lemma 1.6, $N$ is conical and $\bar e$ is strongly
unitary.

Now let $\sim$ be the congruence defined on $N$ by
$$x\sim y\Leftrightarrow(\exists p\in{\bb P})(px=py).$$
It is obvious that $N/\!\!\sim$ is conical.
By construction, $N/\!\!\sim$ is ${\bb P}$-torsion-free.
Let $\pi:N\to N/\!\!\sim$ be the natural projection and let
$j=\pi|_M$. Since $M$ is ${\bb P}$-torsion-free, $j$ is
one-to-one. Since $M$ is cofinal in $N$, $j$ has cofinal image.
Let $x,y\in M$ and $z\in N$ such that $j(x)+\pi(z)=j(y)$. This
means that there exists $p\in{\bb P}$ such that $p(x+z)=py$,
{\it i.e.}, $px+pz=py$. Since $N$ is a unitary extension of $M$,
it follows that $pz\in M$. Therefore, by strong unitarity,
$z\in M$, whence $\pi(z)\in j[M]$: hence, $j$ is unitary.
Finally, let $z\in N$ and $m\in{\bb N}$ such that
$m\pi(y)\in j[M]$. This means that there exists $x\in M$ such
that $my\sim x$, thus there exists $p\in{\bb P}$ such that
$pmy=px$; thus $pmy\in M$, thus, since $e$ is strongly unitary,
$y\in M$, whence $\pi(y)\in j[M]$. Therefore, $j$ is strongly
unitary. Thus identify $M$ with $j[M]$.

But the equation $\Sigma$ obviously admits a
solution in $N$ ({\it viz.} $(x=[(1,0),0];y=[(0,1),0])$
with the notations of Lemma 1.6), thus
also in $N/\!\!\sim$. Since $N/\!\!\sim$ is a
${\bb P}$-torsion-free strong unitary (conical) extension of $M$
and by assumption on $M$, the equation also admits a solution
in $M$. Hence we have checked quasi-divisibility of
$M$.\smallskip

Finally, the fact that $M$ is a refinement cone results from
Theorem 1.14 (a).\qed\medskip

\Corollary3.6. Every ${\bb P}$-torsion-free
cone embeds strongly unitarily into a ${\bb P}$-torsion-free
quasi-divisible refinement cone.\qed\medskip

Again, the proof of Corollary 1.12 allows us to extend this
result to other classes:\medskip

\Corollary3.7. Every
stably finite ${\bb P}$-torsion-free cone embeds strongly
unitarily into a stably finite ${\bb P}$-torsion-free
quasi-divisible refinement cone.\qed\medskip

\Corollary3.8. Every directed
${\bb P}$-torsion-free \poag\ $A$ embeds cofinally into a directed
${\bb P}$-torsion-free Riesz group $B$ with quasi-divisible
positive cone, with in addition $B/A$ torsion-free.
\Proof. Let {\bf C} be the class of all
cancellative ${\bb P}$-torsion-free cones. It suffices to prove
that every {\bf C}-strongly unitarily closed element $M$ of
{\bf C} is a quasi-divisible refinement cone. By Corollary
3.6, $M$ embeds strongly unitarily into a ${\bb P}$-torsion-free
quasi-divisible refinement cone $N$. Now
define a congruence $\equiv$ on $N$ by putting
$$x\equiv y\Leftrightarrow(\exists z)(x+z=y+z).$$

As in the proof of the Claim of Theorem 1.9, one proves that
$N/\!\!\equiv$ is conical.
By using the fact that $M$ is strongly unitary (in particular,
cofinal) in $N$, it is easy to verify that the natural
homomorphism $j:M\to N/\!\!\equiv$ is strongly unitary. Now
identify $M$ and $j[M]$. By using the fact that $N$ is
${\bb P}$-torsion-free, it is also easy to verify that
$N/\!\!\equiv$ is still ${\bb P}$-torsion-free. Thus $N/\!\!\equiv$
belongs to {\bf C}.

The rest of the proof is automatic: let for example $a\in M$.
Since $N$ is quasi-divisible, the equation
$2{\bf x}+3{\bf y}=a$ admits a solution in $N$, thus also in
$N/\!\!\equiv$. Since the latter is a strong unitary extension
of $M$ and by assumption on $M$, the equation also admits a
solution in $M$. The proof for refinement is
similar.\qed\medskip

Let us now turn our attention on {\it simple}
${\bb P}$-torsion-free cones. The ``simple analogue" of Theorem
3.5 is the following:\medskip

\Theorem3.9. Let $M$ be
simply ${\bf C}_{\bb P}$-strongly unitarily closed. Then $M$ is
a quasi-divisible refinement cone.
\Proof. By Corollary 3.6, $M$ embeds
strongly unitarily into some ${\bb P}$-torsion-free
quasi-divisible refinement cone $N$.\smallskip 

Let us first prove quasi-divisibility of $M$.
Thus let $a\in M$, we solve in $M$ the equation
$\Sigma:2{\bf x}+3{\bf y}=a$. If $a=0$ then there is trivially
a solution, so suppose that $a\ne0$. Since $N(a)$ is
quasi-divisible (Lemma 3.4), $\Sigma$ admits a solution in
$N(a)$. Since $a\ne0$ and $M$ is simple, $M\subseteq N(a)$,
thus $N(a)$ is a strong unitary extension of $M$, and since
$N(a)$ is simple and belongs to ${\bf C}_{\bb P}$, $\Sigma$ admits
a solution in $M$. Therefore, $M$ is quasi-divisible.\smallskip

Let us check now that $M$ is a refinement monoid.
Thus let $a_i$, $b_i$ ($i<2$) be elements of $M$ such that
$a_0+a_1=b_0+b_1$; consider as usual the following equation
system
$$\Sigma':\cases{
{\bf x}_{i0}+{\bf x}_{i1}=a_i&(all $i<2$)\cr
{\bf x}_{0i}+{\bf x}_{1i}=b_i&(all $i<2$)\cr
}$$

If $a_i=0$ or $b_i=0$ for some $i$, then it is obvious that
$\Sigma$ admits a solution in $M$. Thus suppose that all the
$a_i$, $b_i$'s are non zero. Then, since $M$ is simple,
all four of them are order-units, thus
$a_0\asymp a_1\asymp b_0\asymp b_1$. By Lemma 3.4, $N(a_0)$ is
a refinement cone, thus $\Sigma'$ admits a solution in
$N(a_0)$. Thus as before, $\Sigma'$ admits a solution in $M$,
and $M$ is a refinement cone.
\qed\medskip

\Corollary3.10. Every
simple ${\bb P}$-torsion-free cone embeds strongly unitarily into
a simple ${\bb P}$-torsion-free quasi-divisible refinement
cone.\qed\medskip

Furthermore, as for Corollary 3.7, one obtains easily\medskip

\Corollary3.11. Every simple
${\bb P}$-torsion-free stably finite cone embeds strongly unitarily
into a simple ${\bb P}$-torsion-free stably finite quasi-divisible
refinement cone.\qed\medskip

Now, a proof similar to the proof of Corollary 3.8 yields
easily the following:\medskip

\Corollary3.12. Every simple
${\bb P}$-torsion-free \poag\ $A$ embeds cofinally into a simple
${\bb P}$-torsion-free Riesz group with quasi-divisible positive
cone $B$ such that in addition, $B/A$ is
torsion-free.\qed\medskip

\Remark3.13. In fact, a direct (though not
devoid of lengthy calculations...) construction shows that
{\it every \poag\ (resp. simple \poag) $G$ embeds cofinally
into an interpolation group (resp. simple Riesz group) $H$ such
that in addition, if $F$ is the free abelian group with $|G^+|$
generators, then $H=G\oplus F$ as abelian groups}. Note that
this implies immediately that if $G$ is torsion-free, then so
is $H$.\medskip

Corollary 3.12 also allows us to construct the following
example (the question which it answers was communicated to us
by K. R. Goodearl):\medskip

\Example3.14. {\sl A countable torsion-free
simple Riesz group $G$ and an interval ({\rm i.e.}, a nonempty,
upward directed lower subset of $G^+$) $D$ of $G^+$ such that
$2D=G^+$ but $D\ne G^+$.}
\Proof. We shall first construct an example with
all the properties above except interpolation; then we will
conclude by Corollary 3.12.\smallskip

To start with, let $A$ be the submonoid of ${\bb Z}^+$ generated
by $\{2,7\}$: that is,
$$A=\{0,2,4,6,7,8,9,10,11,\ldots\}$$
and let $M$ be the submonoid of ${\bb Q}^+$ generated by all
elements of the form $(k/2)(9/2)^n$ where $k\in A$ and
$n\in{\bb Z}^+$. Let $G_0=M+(-M)$, equipped with the positive cone
$M$. Since $G_0$ is directed and $G_0^+=M$ is a submonoid of
${\bb Q}^+$, {\it$G_0$ is a simple \poag}.
For all $n\in{\bb Z}^+$, put $d_n=(9/2)^n$.
Note that $d_n=(2/2)(9/2)^n$, 
$d_{n+1}-d_n=(7/2)(9/2)^n$ and
$2d_{n+1}-4d_n=(10/2)(9/2)^n$, whence the following claim:
\medskip

\noindent{\bf Claim 1.} {\sl For all $n\in{\bb Z}^+$, all elements
$d_n$, $d_{n+1}-d_n$ and $2d_{n+1}-4d_n$ belong to
$M$.}\qed\medskip

Since, by Claim 1, the sequence $(d_n)_n$ is increasing in
$G_0$, it generates an interval $D_0$, {\it viz.}
$$D_0=\{x\in G_0^+:\ (\exists n\in{\bb Z}^+)(x\leq_{G_0}d_n)\}.$$
\medskip

\noindent{\bf Claim 2.} $2D_0=G_0^+$.\smallskip

\noindent{\bf Proof of Claim.} It is easy to see that every
element of $G_0$ is bounded above (for $\leq_{G_0}$) by
some $m(d_0+d_1+\cdots+d_n)$ where $m,n\in{\bb N}$, thus by
$mnd_n$, thus by some $2^kd_n$ ($k\in{\bb N}$). But an easy
induction proof (using Claim 1) shows that
$2^kd_n\leq_{G_0}2d_{n+k-1}$; since $2d_{n+k-1}\in2D_0$, the
conclusion follows.\qed\medskip

We shall now prove that $2d_0\notin D_0$ (whence
$D_0\ne G_0^+$). Towards this goal, we shall prove by
induction on $m$ that $2d_0\not\leq_{G_0}d_m$ (all
$m\in{\bb Z}^+$). It is trivial for $m=0$. Suppose that $m>0$
and that the result has been proved for all $m'<m$, and
suppose that the conclusion fails for $m$, {\it i.e.}, there are
$n\in{\bb Z}^+$ and $k_l$ ($l\leq n$) in $A$ such that
$$(9/2)^m-2=\sum_{l\leq n}(k_l/2)(9/2)^l.\leqno(*)$$
Taking the minimal possible value for $n$ ensures that
$k_n\ne0$ (because $(9/2)^m\ne2$). Then
$(9/2)^m-2\geq(2/2)(9/2)^n$, whence $n<m$. On the other hand,
the right-hand side of $(*)$ belongs to $2^{-n-1}{\bb Z}^+$, thus
$2^{n+1}((9/2)^m-2)\in{\bb Z}^+$, thus $n+1\geq m$. It follows that
$m=n+1$. If $k_n$ were even, then the right-hand side of $(*)$
would belong to $2^{-n}{\bb Z}^+$, whence
$2^n((9/2)^{n+1}-2)\in{\bb Z}^+$, a contradiction: thus $k_n$ is
odd. If $k_n\geq9$, then the right-hand side of $(*)$ would be
$\geq(9/2)^{n+1}$, thus $>(9/2)^m-2$, a contradiction: thus the
only possibility left is $k_n=7$, so that
$(k_n/2)(9/2)^n=(9/2)^m-(9/2)^n$. Hence, after canceling
$(9/2)^m$ from $(*)$, we obtain
$$(9/2)^n-2=\sum_{l<n}(k_l/2)(9/2)^l,$$ with $n<m$, which
contradicts the induction hypothesis. Thus we have proved that
$2d_0\notin D_0$.\smallskip

By Corollary 3.12, $G_0$ embeds cofinally into a
torsion-free simple Riesz group $G$; a
standard L\"owenheim-Skolem type argument shows easily that
one can take $G$ {\it countable}. If $D$ is the interval of
$G^+$ generated by $D_0$, it is then easy to see that $2D=G^+$
but $D\ne G^+$.\qed\medskip

\noindent{\bf Problem 3.15.} Can one realize Example 3.14 as
a torsion-free Riesz group of rank one ({\it i.e.}, with
positive cone an additive submonoid of ${\bb Q}^+$)?
The latter were for example studied in [{\bf\pard}].\medskip

\noindent{\bf Final comments.} As mentioned in the
Introduction, we chose here to restrict ourselves essentially
to {\it conical} \cm s and all our results here can be
extended, modulo sometimes minor changes in the statements and
the proofs, to arbitrary \cm s (for example, Lemma 1.6 (a) is
no longer useful but some extra care has to be taken about the
maximal subgroup of our monoids). The definition of $M(a)$ has
then to be changed into
$M(a)=\{x\in M:\ x\asymp a\ \or\ x\leq0\}$.
Some results whose
statements extend {\it mutatis mutandis} are for example
1.8, 1.9, 1.10, 1.12, 1.14, 1.15 (a)(c), 2.6, 2.7, 2.8, 2.10,
3.5 -- 3.12. The statements of 1.16, 1.17, 2.12 also extend to
all \cm s, but antisymmetry is no longer needed among the
hypotheses.\medskip

\noindent{\bf Acknowledgment.} I wish to thank Ken Goodearl
for his many helpful comments and corrections about this
paper, and in particular for having pointed several oversights
and a significant error in its first version.

\references

\lastpage
\bye